# Equivalent quasi-norms on rearrangement-invariant quasi-Banach function spaces

Leo R. Ya. Doktorski[1*]

[1*]Fraunhofer IOSB (Institute of Optronics, System Technologies and Image Exploitation), Department Object Recognition, Gutleuthausstr. 1, 76275, Ettlingen, Germany.
*Corresponding author. E-mails: doktorskileo@gmail.com; leo.doktorski@iosb.fraunhofer.de

## Abstract

We present new formulae providing equivalent quasi-norms on rearrangement-invariant quasi-Banach function spaces. Our results are based on properties of certain averaging operators and formulated in terms of Boyd indices. We present an illustrative application of our general result to the scales of Lorentz–Karamata spaces, so-called $\mathcal{L}$- and $\mathcal{R}$-Lorentz–Karamata spaces and to grand and small Lorentz spaces.

## 1 Introduction

Let $(\Omega, \mu)$ be a totally $\sigma$-finite measure space with a non-atomic measure $\mu$. By $\mathfrak{M}(\Omega, \mu)$ we denote the set of all $\mu$-measurable functions on $\Omega$. If $f \in \mathfrak{M}(\Omega, \mu)$, $f^*(t)$ ($t$>0) denotes the non-increasing rearrangement of $f$ and

$$f^{**}(t) := \frac{1}{t}\int_0^t f^*(u)du \tag{1}$$

is the maximal function of $f^*$. (See, e.g., [5].) The most common way to define a rearrangement-invariant (abbreviation r.i.) quasi-Banach function space on $\Omega$ is to define an expression for the quasi-norm containing $f^*$ or $f^{**}$. The machinery goes back to Lorentz [53] and is used in most papers and books dealing with r.i. quasi-Banach function spaces. See, for example, [5, 6, 51, 69] and references therein. The same expression written with $f^*$ and with $f^{**}$ can produce two quasi-norms that define, a priori, different quasi-Banach function spaces. Thus, it is natural to ask: How are the resulting spaces related to one another and when do they coincide?

For example, if $0 < p, q \le \infty$, the Lorentz spaces $L_{p,q}$ and $L_{(p,q)}$ are defined as the sets of all $f \in \mathfrak{M}(\Omega, \mu)$ such that

$$\|f\|_{p,q} := \left\| t^{\frac{1}{p}-\frac{1}{q}} f^*(t) \right\|_{q,(0,\infty)} < \infty$$

or

$$\|f\|_{(p,q)} := \left\| t^{\frac{1}{p}-\frac{1}{q}} f^{**}(t) \right\|_{q,(0,\infty)} < \infty,$$

respectively. (See, for example, [5, Chapter 4, Definitions 4.1 and 4.4])
It is known that if $1 < p \le \infty$ and $0 < q \le \infty$, then the quasi-norms $\|*\|_{p,q}$ and $\|f\|_{(p,q)}$ are equivalent. (See, e.g., [5, Chapter 4, Lemma 4.5] and [64, Theorems 3.8 and 3.16]). Similar equivalences are obtained for some other scales of spaces. See, e.g., [23, 64, 65] and references therein. The establishing of such equivalences for many scales of function spaces arising in the $K$-interpolation method in limiting cases is also the main aim of the paper [23].

Formula (1) can be considered as an averaging (Hardy) operator. Thus, it is possible to formulate a more general problem of searching for other equivalent quasi-norms on r.i. quasi-Banach function spaces using various averaging operators.

In [29], new characterizations of Lorentz–Karamata spaces were given by means of quasi-norms that were shown to be equivalent to the classical ones. The results are based on properties of certain averaging operators. Earlier, similar results for Lorentz spaces are presented in [28, 62, 63]. The technique used here is based on Hardy-type inequalities. However, Lorentz–Karamata spaces (and their special cases, Lorentz spaces and Lorentz–Zygmund spaces) are only one scale in a family of so-called r.i. quasi-Banach function spaces that have proved useful in analysis.

The aim of this paper is to establish new formulae providing equivalent quasi-norms on arbitrary r.i. quasi-Banach function spaces. The results are also based on properties of certain averaging operators. The necessary (and sometimes the sufficient) conditions of equivalence for different quasi-norms are formulated in terms of Boyd indices.

The computation of Boyd indices for various r.i. quasi-Banach function spaces is a topic of many papers. See, e.g., [9, 40, 48, 58, 65] and references therein. In our applications, we calculate Boyd indices using known results related to "simpler" spaces and some easy interpolation arguments.

The organization of the paper is the following. Section 2 contains the basic tools: definitions and some properties of slowly varying functions, quasi-normed spaces, rearrangement-invariant quasi-Banach function spaces and Boyd indices. Subsection 2.6 is devoted to averaging (Hardy-type) operators we will be using. Main results of this paper dealing with equivalent quasi-norms on r.i. quasi-Banach function spaces are proved in Section 3. In Section 4, we recall some basic definitions of interpolation functors and real interpolation spaces we shall work with. As applications (see Section 5.) we reformulate the general results, formulated in terms of Boyd indices, in terms of parameters of concrete functional scales. In Subsection 5.1 we consider Lorentz–Karamata spaces. The main results of this subsection improve slightly main results of [29]. In Subsections 5.2 and 5.3 we study so-called $\mathcal{L}$- and $\mathcal{R}$-Lorentz–Karamata spaces and grand and small Lorentz spaces involving slowly varying functions, respectively.

## 2 Preliminaries

For $f$ and $g$ being non-negative functions, we write $f \prec g$ if $f \leq Cg$, where the positive constant $C$ is independent on all significant quantities. Two non-negative functions $f$ and $g$ are considered equivalent ($f \sim g$) if $f \prec g$ and $g \prec f$.

Let $X$ and $Y$ be two quasi-normed spaces with quasi-norms $\|*\|_X$ and $\|*\|_Y$, respectively. Throughout the paper, we write $X \hookrightarrow Y$ to indicate that $X$ is continuously embedded in $Y$, i.e. there exists a constant $C > 0$ such that $\|x\|_Y \leq C\|x\|_X$ for all $x$ in $X$. We write $X = Y$ if $X \hookrightarrow Y$ and $Y \hookrightarrow X$. In this case, $X$ and $Y$ are equal as sets and as linear spaces and their quasi-norms are equivalent ($\|x\|_X \sim \|x\|_Y$). We write $T: X \to Y$ to indicate that $T$ is a bounded operator from $X$ into $Y$. Its norm we denote by $\|T\|_{X\to Y}$.

Throughout the paper, $L \in \{1, \infty\}$. By $\mathfrak{M} \equiv \mathfrak{M}(0, L)$ we will denote the set of all complex-valued, Lebesgue measurable and finite almost everywhere (abbreviation: a.e.) functions on $(0, L)$. As it is customary, we will identify functions that coincide a.e.. We will further denote by $\mathfrak{M}^+$ the subset of $\mathfrak{M}$ containing the non-negative functions. Given a measurable set $\omega \subset (0, L)$ we will denote its characteristic function by $\chi_\omega$.

Although we are considering functions on $(0, L)$, because of the Luxemburg representation theorem all our results can be reformulated for $\mu$-measurable functions defined on a $\sigma$-finite measure space $(\Omega, \mu)$ with a non-atomic measure $\mu$. See [56] and [60, Theorem 3.1 and Proposition 3.3] for r.i. quasi-Banach function spaces. Note that these measure spaces are resonant. See [5, Chapter 2, Definition 2.3, p.45 and Theorem 2.7, p 52].

Let $f$ be a real-valued non-negative function on $(0, L)$. We say that it is *almost* increasing (decreasing) if it is equivalent to a non-decreasing (non-increasing) function.

We adopt the conventions $1/\infty = 0$ and $1/0 = \infty$.

A few more pieces of notation:

$\|*\|_{p,(a,b)}$ is the usual (quasi-)norm on the classical Lebesgue space $L_p$ on the interval $(a, b)$ ($0 < p \leq \infty$, $0 \leq a < b \leq \infty$). If $a = 0$ and $b = \infty$, we write $\|*\|_p$.

$f_n \nearrow f$: the sequence $\{f_n\}$ $(n = 1,2, \ldots)$ of real-valued non-negative functions is pointwise monotone increasing and we have the pointwise convergence $f_n \to f$ a. e..

## 2.1 Slowly varying functions

In this subsection, we summarize some properties of so-called slowly varying functions, which will be required later. For more details, we refer to, e.g., [26, 32, 34, 44, 65, 66].

**Definition 1** We say that a positive Lebesgue measurable function $b$ on $(0,L)$ is slowly varying, notation $b \in SV$, if, for each $\varepsilon > 0$, the function $t^\varepsilon b(t)$ is an almost increasing function while the function $t^{-\varepsilon} b(t)$ is an almost decreasing function.

Examples of $SV$–functions include functions which are equivalent to 1, powers of logarithms

$$\ell^\alpha(t) = (1 + |log t|)^\alpha, \qquad \alpha \in \mathbb{R},$$

reiterated logarithms

$$(\ell \circ \ldots \circ \ell)^\alpha(t), \qquad \alpha \in \mathbb{R},$$

"broken logarithmic" functions of [30, 31, 64] and the family of functions $\exp(|log t|^\alpha)$, for $\alpha \in (0,1)$. Some basic properties of slowly varying functions are summarized in the following lemma. Throughout the paper, we will often use them without any explicit reference.

**Lemma 2** *Let $b, b_1, b_2 \in SV$, $\alpha > 0$, $0 < q \leq \infty$, $\lambda \in (-\infty, \infty)$ and $t \in (0, L)$.*

(i) *Then $b^\lambda \in SV$, $b\left(\frac{1}{t}\right) \in SV$, $b\big(t^\alpha b_1(t)\big) \in SV$ and $b_1 b_2 \in SV$.*

(ii) *If $f \sim g$, then $b \circ f \sim b \circ g$.*

(iii) $\left\|u^{\alpha - 1/q} b(u)\right\|_{q,(0,t)} \sim t^\alpha b(t)$ *and* $\left\|u^{-\alpha - 1/q} b(u)\right\|_{q,(t,\infty)} \sim t^{-\alpha} b(t)$.

(iv) *The functions $\left\|u^{-1/q} b(u)\right\|_{q,(0,t)}$ and $\left\|u^{-1/q} b(u)\right\|_{q,(t,\infty)}$ (if they exist) belong to SV and*

$b(t) \prec \left\|u^{-1/q} b(u)\right\|_{q,(0,t)}$ *and* $b(t) \prec \left\|u^{-1/q} b(u)\right\|_{q,(t,\infty)}$.

(v) $\left\|u^{\lambda - 1/q} b(u)\right\|_{q,\left(\frac{t}{2}, t\right)} \sim \left\|u^{\lambda - 1/q} b(u)\right\|_{q,(t,2t)} \sim t^\lambda b(t)$.

It follows from property (iii) that every slowly varying function is equivalent to a continuous one. Furthermore, it is also equivalent to an infinitely differentiable function, see [66].

### 2.2 Quasi-normed spaces

We recall some basic facts concerning quasi-normed spaces. As reference for this subject, we recommend [6].

**Definition 3** A quasi-norm $\|x\|_X$ defined on a vector space $X$ (over a field $\mathbb{K}$ of real or complex numbers) is a map $X \to \mathbb{R}^+$ such that

(i) $\|x\|_X = 0 \iff x = 0$,

(ii) $\|\lambda x\|_X = |\lambda| \|x\|_X$ for $\lambda \in \mathbb{K}$, $x \in X$,

(iii) $\|x + y\|_X \leq C_X(\|x\|_X + \|y\|_X)$ for all $x, y \in X$, where $C_X$ is a constant independent of $x, y$.

A vector space $X$ with a quasi-norm $\|x\|_X$ is called a quasi-normed space.

If it is important to emphasize that we are considering the space $X$ with the quasi-norm $\|*\|_X$, we will write $\{X, \|*\|_X\}$. If $X$ is complete, it will be called quasi-Banach space. The inequality (iii) above is called the "*quasi-triangle inequality*" and the constant $C_X$ is called "*quasi-triangle constant*" or "*the modulus of concavity*". Obviously $C_X \geq 1$. If $C_X = 1$, then $X$ is a normed space.

**Lemma 4** ([6, Lemma 3.10.3]) *The classical Lebesgue space* $L_p$ *is a Banach space if* $1 \leq p \leq \infty$ *and a quasi-Banach space if* $0 < p < 1$. *It holds*

$$\|f + g\|_{L_p} \leq \max\left\{1, 2^{\frac{1-p}{p}}\right\}\left(\|f\|_{L_p} + \|g\|_{L_p}\right).$$

### 2.3 Quasi-Banach function spaces

Banach function spaces were introduced in 1955 by Luxemburg [55]. We need them as an intermediate step in defining r.i quasi-Banach function spaces in the next subsection. Despite the vast number of works devoted to (quasi-)Banach functional spaces, there are various definitions of this concept in the literature. See, for example, [9, 55, 56, 58, 61], [5, Chapter 1, Definition 1.1], [11, Definition 3.1] and [60, Definition 2.10.] and references therein. A discussion of various definitions of quasi-Banach functional spaces can be found in [54]. A little different approach via so-called ideal spaces can be found, e.g., in [49, 50, 52]. Our definition differs slightly from the well-known definitions. See *Remark 6* below.

**Definition 5** A mapping $\rho: \mathfrak{M}^+(0, L) \to [0, \infty]$ ($L \in \{1, \infty\}$) is called a function quasi-norm if there exists a constant $C_\rho \geq 1$ such that for $f, g, f_n (n \in \mathbb{N})$ in $\mathfrak{M}^+(0, L)$ and for all constants $\lambda \geq 0$ the following properties hold:

(P1) $\rho(f) = 0 \Leftrightarrow f = 0$ a.e.; $\rho(\lambda f) = \lambda\rho(f)$; $\rho(f + g) \leq C_\rho\big(\rho(f) + \rho(g)\big)$;

(P2) $\rho$ has the lattice property or ideal property: if $g \leq f$ a.e., then $\rho(g) \leq \rho(f)$;

(P3) $\rho$ has the Fatou property: if $0 \leq f_n \nearrow f$ a.e. and $\sup_{n\geq 1} \rho(f_n) < \infty$, then $\rho(f) < \infty$ and $\rho(f_n) \nearrow \rho(f)$;

(P4-weak) If a set $\omega \subset (0, L)$ is bounded and measurable, then $\rho(\chi_\omega) < \infty$.

The corresponding *quasi-Banach function space* $\{B_\rho, \rho(|f|)\}$ we define as the set

$$B_\rho \equiv B_\rho(0, L) = \left\{f \in \mathfrak{M}(0, L): \ \|f\|_{B_\rho} := \rho(|f|) < \infty\right\}.$$

We will write $B_\rho(0, L)$ only if we want to emphasize that the functions we are considering are defined on $(0, L)$. It is known that the quasi-Banach function spaces are complete. See [11, Lemma 3.6], [61, Corollary 3.8] and [54, Proposition 2.2]. Consequently, it is correct to refer to these spaces as "Banach spaces."

*Remark 6* Our definition of quasi-Banach function space differs slightly from the standard one. Cf., e.g., [5, 8, 21, 55, 59, 60, 61]. In the literature, instead of the property (P4-weak) is most often formulated the following one.

(P4) If the set $\omega \subset (0, L)$ is measurable and $\operatorname{mes} \omega < \infty$, then $\rho(\chi_\omega) < \infty$.

Obviously, (P4-weak) and (P4) coincide if $L = 1$. However, it may be too restrictive if $L = \infty$. A discussion of this problem can be found in [68]. Example 10 below illustrates our point. Note that if one considered function spaces on $\mathbb{R}^n$, then the property (P4-weak) coincides with the property (iv) in [68, Definition 2.1] where the ball quasi-Banach function spaces are defined and studied.

## 2.4 Rearrangement-invariant quasi-Banach function spaces

In this subsection we recall the definition of rearrangement-invariant quasi-Banach function spaces. R.i. spaces are spaces whose norms depend only on the non-increasing rearrangement of measurable functions. They include not only the Lebesgue spaces mentioned above but also the Lorentz, Zygmund, Lorentz–Zygmund, generalized Lorentz–Zygmund, Lorentz–Karamata and Orlicz spaces, as well as others.

Recall that throughout this paper we consider complex-valued, Lebesgue measurable and finite a.e. functions on $(0, L)$, $L \in \{1, \infty\}$. The *non-increasing rearrangement* $f^*$ of $f \in \mathfrak{M}$ is defined by (see, e.g., [5])

$$f^*(t) = \{\lambda > 0: \mu\{x \in (0, \infty): |f(x)| > \lambda\} \leq t\}, \qquad t > 0.$$

It is known that the non-increasing rearrangement is unique. The functions $f$ and $g$ are called equimeasurable if $f^* = g^*$. More information about this topic can be found, e.g., in [5, Chapter 2].

**Definition 7** Let $\rho: \mathfrak{M}^+(0, L) \to \mathbb{R}^+$ be a function quasi-norm and $B_\rho$ be the corresponding quasi-Banach function space. The set $E_\rho$ defined as

$$E_\rho \equiv E_\rho(0, L) = \left\{ f \in \mathfrak{M}: \ \|f\|_{E_\rho} \coloneqq \rho(f^*) = \|f^*\|_{B_\rho} < \infty \right\}$$

is called *rearrangement-invariant* (abbreviation *r.i.*) *quasi-Banach function space*. Sometimes we will denote it as $\{E_\rho, \rho(f^*)\}$.

*Remark 8* Note that for the space $E_\rho$, the property (P4) holds. Really, let $\omega \subset (0, L)$ be a measurable set with finite measure. By (P4-weak) for the space $B_\rho$, we have

$$\|\chi_\omega\|_E = \|(\chi_\omega)^*\|_B = \left\|\chi_{(0,\mu(\omega))}\right\|_B < \infty.$$

*Remark 9* Sometimes in the literature, $E_\rho$ is referred to as the symmetrization of $B_\rho$. See, e.g., [47, 49, 50].

**Definition 7** of r.i. (quasi-) Banach function spaces includes the Fatou property (P3). Many authors (See, for example, [12, 13, 45, 52].) investigate a wider class of spaces of so-called *symmetric (quasi-) Banach function spaces*. The definition of these spaces does not include the Fatou property. The classical monographs [51] and [5] are devoted to symmetric and to r.i. Banach function spaces, respectively. We limit ourselves to considering only r.i. (quasi-) Banach function spaces.

*Example 10* Let $0 < p < r \leq \infty$ and the mapping $\rho_{p,r}: \mathfrak{M}^+(0, \infty) \to [0, \infty]$ be defined by the formula

$$\rho_{p,r}(f) = \left\| t^{\frac{1}{p} - \frac{1}{r}} f(t) \right\|_{r,(0,\infty)} < \infty.$$

Clearly, properties (P1), (P2) and (P3) are fulfilled. Let us check (P4-weak). Consider the case $r < \infty$. The case $r = \infty$ can be considered similarly. If $\omega \subset (0, A)$ where $0 < A < \infty$ is a measurable set, then because $\frac{1}{p} > 0$,

$$\rho_{p,r}(\chi_\omega) = \left\| t^{\frac{1}{p}-\frac{1}{r}} \chi_\omega(t) \right\|_{r,(0,\infty)} \le \left\| t^{\frac{1}{p}-\frac{1}{r}} \right\|_{r,(0,A)} = \left(\frac{p}{r}\right)^{\frac{1}{r}} A^{\frac{1}{p}} < \infty.$$

So, property (P4-weak) is fulfilled. Let now $\alpha > 1$, $\omega_k = (k^\alpha, k^\alpha + \frac{1}{k^2})$, $(k = 1,2, \ldots)$ and $\omega = \bigcup_{k \ge 1} \omega_k$. Because $\alpha > 1$, the sets $\omega_k$ do not intersect. The set $\omega$ is unbounded and has finite measure:

$$\operatorname{mes} \omega = \sum_{k \ge 1} \frac{1}{k^2} < \infty.$$

Moreover, because $\frac{r}{p} > 1$,

$$\{\rho_{p,r}(\chi_\omega)\}^r = \sum_{k \ge 1} \int_{k^\alpha}^{k^\alpha + \frac{1}{k^2}} t^{\frac{r}{p}-1} dt \ge \sum_{k \ge 1} k^{\alpha\left(\frac{r}{p}-1\right)} \int_{k^\alpha}^{k^\alpha + \frac{1}{k^2}} dt = \sum_{k \ge 1} k^{\alpha\left(\frac{r}{p}-1\right)-2}.$$

Hence, $\rho_{p,r}(\chi_\omega) = \infty$ if $\alpha$ is large enough. This means that the property (P4) is not fulfilled. Denote by $B_{p,r}$ the Banach function space generated by $\rho_{p,r}$. In this example, it is shown that $\chi_\omega \notin B_{p,r}$ but $\chi_\omega \in L_{p,r}$. (See Definition 29.) Nevertheless, it is clear that $L_{p,r}$ is the symmetrization of $B_{p,r}$.

Let $\rho: \mathfrak{M}^+(0, L) \to \mathbb{R}^+$ be a function quasi-norm, and let $B_\rho$ and $E_\rho$ be the corresponding quasi-Banach function space and r.i. quasi-Banach function space with the quasi-norms

$$\|f\|_{B_\rho} = \rho(|f|) \quad \text{and} \quad \|f\|_{E_\rho} = \rho(f^*) = \|f^*\|_{B_\rho},$$

respectively. Example 10 shows that it can happen that the quasi-norms $\|f\|_{B_\rho}$ and $\|f\|_{E_\rho}$ are not equal and even not equivalent.

## 2.5 Dilation operators and Boyd indices

An important property of r.i. quasi-Banach function spaces is that the dilation operator is bounded on these spaces. This makes it possible to define the Boyd indices. They were introduced by Boyd in [7, 8, 9] and are useful in various problems concerning r. i. spaces. In this subsection, we collect some facts about dilation operators and Boyd indices which will be needed later. To the readers interested in this topic we would recommend monographs [58], [5, Chapter 3, Section 5] and [51, Chapter 2]. New results can be found in [12, 61, 60, 65] and references therein. See also [10, 13, 47, 50, 52, 57, 58]. Recall that $L = 1$ or $L = \infty$.

**Definition 11** Let $a > 0$. The dilation operator $D_a$ on $\mathfrak{M}(0, L)$ is defined by

$$(D_a f)(t) := f(at)\chi_{(0,L)}(at), \qquad 0 < t < L.$$

Obviously,

$$D_a(\lambda f) = \lambda D_a(f), \quad D_a(f + g) = D_a(f) + D_a(g) \quad \text{and} \quad D_a f^* = (D_a f)^*.$$

Let $E$ be some r.i. quasi-Banach function space on $(0, L)$ and $a > 0$. It is known that the dilation operator $D_a$ is bounded on $E$ and the function $a \rightsquigarrow \|D_a\| \equiv \|D_a\|_{E \to E}$ is submultiplicative. This means that

$$\|D_{ab}\| \le \|D_a\| \|D_b\| \qquad \text{for all } 0 < a, b < \infty.$$

Because $f^*$ is non-increasing, $\|D_a\| \le 1$ if $a \ge 1$. The boundedness of dilation operators entails the following lemma.

**Lemma 12** *Let $E$ be a r.i. quasi-Banach function space on $(0, L)$ and $a > 0$. Then for all $f \in E$*

$$\|D_a f\|_E \sim \|f\|_E$$

*with the constants depending on $E$ and $a$.*

We use the definition of the lower and upper Boyd indices in the following form. See, e.g., [12, 13, 52, 59].

**Definition 13** Let $E$ be a r.i. quasi-Banach function space on $(0, L)$. Define the lower Boyd index $p_E$ of $E$ by

$$p_E = \sup\left\{p > 0 : \exists c > 0 \text{ such that } \forall 0 < a < 1 \ \|D_a\| \leq c a^{-\frac{1}{p}}\right\}$$

and the upper Boyd index $q_E$ of $E$ by

$$q_E = \inf\left\{q > 0 : \exists c > 0 \text{ such that } \forall a \geq 1 \ \|D_a\| \leq c a^{-\frac{1}{q}}\right\}.$$

In [12, Proposition 4.6] it is shown that

$$p_E = \sup_{s>1} \frac{\ln s}{\ln\|D_{1/s}\|} = \lim_{s\to\infty} \frac{\ln s}{\ln\|D_{1/s}\|} \quad \text{and} \quad q_E = \sup_{0<s<1} \frac{\ln s}{\ln\|D_{1/s}\|} = \lim_{s\to 0+} \frac{\ln s}{\ln\|D_{1/s}\|}.$$

In many texts (see, e.g., [5, Chapter 3, Definitions 5.10 and 5.12]), the lower and upper Boyd indices of $E$ are alternatively defined as the quantities

$$\underline{\alpha}_E = \sup_{0<s<1} \frac{\ln\|D_{1/s}\|}{\ln s}, \qquad \overline{\alpha}_E = \sup_{s>1} \frac{\ln\|D_{1/s}\|}{\ln s}.$$

It is clear that $\underline{\alpha}_E = \frac{1}{q_E}$ and $\overline{\alpha}_E = \frac{1}{p_E}$. It is known that

$$0 < p_E \leq q_E \leq \infty$$

and if $E$ is a r.i. Banach function space, then

$$1 \leq p_E \leq q_E \leq \infty.$$

## 2.6 Averaging operators

The main results of this paper are based on properties of certain averaging operators and formulated in terms of Boyd indices. Here we define these operators and formulate their properties that we need. We are interested in under what conditions these operators are bounded. The main result of this subsection is **Theorem 19**. It is a corollary of [59, Theorem 2] and [52, Theorem 17].

**Definition 14** (Cf. [14, Definition 5.4].) Let $0 < U, V < \infty$, $0 < W \leq \infty$ and $d \in SV$. For $t \in (0, L)$, we define weighted averaging (Hardy-type) operators on $\mathfrak{M}(0, L)$ as follows:

$$\left(H^{(U,W;d)}f\right)(t) = \left(t^{\frac{1}{U}}d(t)\right)^{-1} \left\|v^{\frac{1}{U}-\frac{1}{W}}d(v)f^*(v)\right\|_{W,(0,t)},$$

$$\left(H_{(V,W;d)}f\right)(t) = \left(t^{\frac{1}{V}}d(t)\right)^{-1} \left\|v^{\frac{1}{V}-\frac{1}{W}}d(v)f^*(v)\right\|_{W,(t,L)},$$

$$\left(H_{(\infty,W;d)}f\right)(t) = \left(d(t)\right)^{-1} \left\|v^{\frac{1}{W}}d(v)f^*(v)\right\|_{W,(t,L)}.$$

These operators will play a crucial role in the future. Without stipulating this each time, we will assume that

$$\begin{cases} d \in SV, \\ 0 < U < \infty, \\ 0 < V \leq \infty, \\ 0 < W \leq \infty. \end{cases} \tag{2}$$

Similar operators with logarithmic term $d$ are used in [14]. If $d = 1$, we skip the character $d$ in the above notations and write, for example, $H^{(U,W)}$. These operators are studied in [59, 52]. In [5, Chapter 3, Definition 5.14] are considered averaging operators

$$(P_\alpha f)(t) = t^{-\alpha} \int_0^t v^\alpha f^*(v) \frac{dv}{v}, \quad (0 < \alpha \le 1)$$

and

$$(Q_\alpha f)(t) = t^{-\alpha} \int_t^\infty v^\alpha f^*(v) \frac{dv}{v}, \quad (0 \le \alpha < 1).$$

Obviously, $P_\alpha = H^{\left(\frac{1}{\alpha},1\right)}$ and $Q_\alpha = H_{\left(\frac{1}{\alpha},1\right)}$. Slightly generalized averaging operators are studied in [57]. For $0 < U < \infty$ put

$$f_{(U)}^{**}(t) := \left(t^{-1} \int_0^t f^*(v)^U dv\right)^{\frac{1}{U}}.$$

These operators were used in [67, 70]. Cf. also [23, (4.2)]. Obviously,

$$f^{**}(t) = f_{(1)}^{**}(t) = \left(H^{(1,1)} f\right)(t) \quad \text{and} \quad f_{(U)}^{**}(t) = \left(H^{(U,U)} f\right)(t).$$

The following lemma describes some properties of the averaging operators from **Definition 14**.

**Lemma 15** (Cf. [70, (R8)] and [23, Lemma 4.1].) *Let* $H$ *be one of the averaging operators* $H^{(U,W;d)}$ *or* $H_{(V,W;d)}$ *from Definition 14 and* (2) *holds. Let* $\lambda \in \mathbb{C}$ *and* $f, g \in \mathfrak{M}(0, L)$. *Then*

(i) $Hf = 0 \Leftrightarrow f = 0$;
(ii) $H(\lambda f) = |\lambda| Hf$;
(iii) $H(|f|) = H(f)$;
(iv) *If* $0 \le |g| \le |f|$ *a.e., then* $Hg \le Hf$;
(v) $H(f^*) = H(f)$.

*Additionally, for all* $t \in (0, L)$ *it holds*

(vi) $\big(H(f+g)\big)(t) \le C_d \max\left\{1, 2^{\frac{1-W}{W}}\right\} \left\{(Hf)\left(\frac{t}{2}\right) + (Hg)\left(\frac{t}{2}\right)\right\}$,

*where the constant* $C_d$ *depends only on* $d$. $C_d = 1$ *if* $d$ *is a constant*;

(vii) $f^*(t) \le \left(H^{(U,\infty;d)} f\right)(t) \prec \left(H^{(U,W;d)} f\right)(t)$, $(0 < W < \infty)$ (Cf. [14, Lemma 5,6].);

(viii) $f^*(2t) \prec \left(H_{(V,W;d)} f\right)(t)$ (Cf. [15, Lemma 7.1].) *and* $\left(H_{(V,\infty;d)} f\right)(2t) \prec \left(H_{(V,W;d)} f\right)(t)$ if $W < \infty$, $f^*(t) \le \left(H_{(V,\infty;d)} f\right)(t)$.

*Proof* The assertions (i-iii) follow from Definition 14. (iv) follows also from Definition 14 because if $0 \le |g| \le |f|$ a.e., then $g^* \le f^*$. (See, e.g., [5, Chapter 2, (1.14)].) (v) follows from the observation that $(f^*)^* = f^*$.

Proof of the assertion (vi). We consider only the case $H = H^{(U,W;d)}$. The case $H = H_{(V,W;d)}$ can be considered analogously. Because $(f+g)^*(s) \le f^*\left(\frac{s}{2}\right) + g^*\left(\frac{s}{2}\right)$ and using a change of variables, we have

$$\begin{aligned}
\left(H^{(U,W;d)}(f+g)\right)(t) &= \left(t^{\frac{1}{U}} d(t)\right)^{-1} \left\| s^{\frac{1}{U}-\frac{1}{W}} d(s) (f+g)^*(s) \right\|_{W,(0,t)} \\
&\le \left(t^{\frac{1}{U}} d(t)\right)^{-1} \left\| s^{\frac{1}{U}-\frac{1}{W}} d(s) \left(f^*\left(\frac{S}{2}\right) + g^*\left(\frac{S}{2}\right)\right) \right\|_{W,(0,t)} \\
&= 2^{\frac{1}{U}} \left(t^{\frac{1}{U}} d(t)\right)^{-1} \left\| v^{\frac{1}{U}-\frac{1}{W}} d(2v) \big(f^*(v) + g^*(v)\big) \right\|_{W,\left(0,\frac{t}{2}\right)}.
\end{aligned}$$

By Lemma 2 (ii), $\exists C_d \geq 1$: $d(2v) \leq C_d d(v)$ and $\frac{1}{d(t)} \leq \frac{C_d}{d\left(\frac{t}{2}\right)}$. Hence, using Lemma 4, we continue

$$\left(H^{(U,W;d)}(f+g)\right)(t) \leq C_d{}^2 \frac{1}{\left(\frac{t}{2}\right)^{\frac{1}{U}} d\left(\frac{t}{2}\right)} \left\| v^{\frac{1}{U}-\frac{1}{W}} d(v)\left(f^*(v)+g^*(v)\right)\right\|_{W,\left(0,\frac{t}{2}\right)}$$

$$\leq C_d{}^2 \max\left\{1, 2^{\frac{1-W}{W}}\right\} \frac{1}{\left(\frac{t}{2}\right)^{\frac{1}{U}} d\left(\frac{t}{2}\right)} \left\{\left\| v^{\frac{1}{U}-\frac{1}{W}} d(v) f^*(v)\right\|_{W,\left(0,\frac{t}{2}\right)} + \left\| v^{\frac{1}{U}-\frac{1}{W}} d(v) g^*(v)\right\|_{W,\left(0,\frac{t}{2}\right)}\right\}$$

$$= C_d{}^2 \max\left\{1, 2^{\frac{1-W}{W}}\right\}\left\{\left(H^{(U,W;d)} f\right)\left(\frac{t}{2}\right) + \left(H^{(U,W;d)} g\right)\left(\frac{t}{2}\right)\right\}.$$

Proof of the assertion (vii). Because $f^*$ is non-increasing, using Lemma 2 (iii), we have

$$f^*(t) \leq f^*(t)\left(t^{\frac{1}{U}} d(t)\right)^{-1} \sup_{0<v<t} v^{\frac{1}{U}} d(v)$$

$$\leq \left(t^{\frac{1}{U}} d(t)\right)^{-1} \sup_{0<v<t} v^{\frac{1}{U}} d(v) f^*(v) = \left(H^{(U,\infty;d)} f\right)(t)$$

$$\sim \left(t^{\frac{1}{U}} d(t)\right)^{-1} \sup_{0<v<t} f^*(v)\left\{\int_0^v \left(s^{\frac{1}{U}} d(s)\right)^W \frac{ds}{s}\right\}^{\frac{1}{W}}$$

$$\leq \left(t^{\frac{1}{U}} d(t)\right)^{-1} \sup_{0<v<t}\left\{\int_0^v \left(s^{\frac{1}{U}} d(s) f^*(s)\right)^W \frac{ds}{s}\right\}^{\frac{1}{W}}$$

$$= \left(t^{\frac{1}{U}} d(t)\right)^{-1}\left\{\int_0^t \left(s^{\frac{1}{U}} d(s) f^*(s)\right)^W \frac{ds}{s}\right\}^{\frac{1}{W}} = \left(H^{(U,W;d)} f\right)(t).$$

Let us prove (viii). We use that $f^*$ is non-increasing and $f^*(v)=0$ if $v>L$. By Lemma 2 (v), we get

$$\left(H_{(V,W;d)} f\right)(t) \geq \left(t^{\frac{1}{V}} d(t)\right)^{-1} \left\| v^{\frac{1}{V}-\frac{1}{W}} d(v) f^*(v)\right\|_{W,(t,2t)}$$

$$\geq \left(t^{\frac{1}{V}} d(t)\right)^{-1} f^*(2t) \left\| v^{\frac{1}{V}-\frac{1}{W}} d(v)\right\|_{W,(t,2t)} \sim \left(t^{\frac{1}{V}} d(t)\right)^{-1} f^*(2t) t^{\frac{1}{V}} d(t) = f^*(2t).$$

Similarly,

$$\left(H_{(V,W;d)} f\right)(t) = \left(t^{\frac{1}{V}} d(t)\right)^{-1} \sup_{2t<v<\infty} \left\| u^{\frac{1}{V}-\frac{1}{W}} d(u) f^*(u)\right\|_{W,\left(\frac{v}{2},\infty\right)}$$

$$\geq \left(t^{\frac{1}{V}} d(t)\right)^{-1} \sup_{2t<v<\infty} \left\| u^{\frac{1}{V}-\frac{1}{W}} d(u) f^*(u)\right\|_{W,\left(\frac{v}{2},v\right)}$$

$$\geq \left(t^{\frac{1}{V}} d(t)\right)^{-1} \sup_{2t<v<\infty} f^*(v) \left\| u^{\frac{1}{V}-\frac{1}{W}} d(u)\right\|_{W,\left(\frac{v}{2},v\right)}$$

$$\sim \left((2t)^{\frac{1}{V}} d(2t)\right)^{-1} \sup_{2t<v<\infty} f^*(v) v^{\frac{1}{V}} d(v) = \left(H_{(V,\infty;d)} f\right)(2t).$$

Further,

$$\left(H_{(V,\infty;d)} f\right)(t) = \left(t^{\frac{1}{V}} d(t)\right)^{-1} \left\| v^{\frac{1}{V}} d(v) f^*(v)\right\|_{\infty,(t,\infty)} \geq \left(t^{\frac{1}{V}} d(t)\right)^{-1} t^{\frac{1}{V}} d(t) f^*(t) = f^*(t). \blacksquare$$

Now, we are interested in the behavior of the functions $(Hf)(t)$ for the averaging operators $H$ from Definition 14 in the case $d=1$.

**Lemma 16** *Let $H$ be one of the averaging operators $H^{(U,W)}$ or $H_{(V,W)}$. Then the function $Hf$ (if it is bounded) is non-increasing and hence $(Hf)^*=Hf$.*
*Proof* This statement is trivial for the operators $H_{(V,W)}$. Consider the operators $H^{(U,W)}$ for $0<U,W<\infty$. Let $\tau=\frac{t}{A}$ with $A>1$. Then, because $f^*$ is non-increasing,

$$\left(H^{(U,W)}f\right)(t)=t^{-\frac{1}{U}}\left\{\int_0^t\left(v^{\frac{1}{U}}f^*(v)\right)^W\frac{dv}{v}\right\}^{\frac{1}{W}}=A^{-\frac{1}{U}}\tau^{-\frac{1}{U}}\left\{\int_0^{A\tau}\left(v^{\frac{1}{U}}f^*(v)\right)^W\frac{dv}{v}\right\}^{\frac{1}{W}}$$

$$=\tau^{-\frac{1}{U}}\left\{\int_0^{\tau}\left(u^{\frac{1}{U}}f^*(Au)\right)^W\frac{du}{u}\right\}^{\frac{1}{W}}\le\tau^{-\frac{1}{U}}\left\{\int_0^{\tau}\left(u^{\frac{1}{U}}f^*(u)\right)^W\frac{du}{u}\right\}^{\frac{1}{W}}=\left(H^{(U,W)}f\right)(\tau).$$

Hence, $\left(H^{(U,W)}f\right)(t)$ is non-increasing. The proof for $H^{(U,\infty)}$ $(0<U<\infty)$ is similar. ■

### 2.7 Boundedness of the averaging operators

Next in this subsection, we assume that $\rho$ is a function quasi-norm on $(0,L)$ and $\{B,\|f\|_B=\rho(|f|)\}$ and $\{E,\|f\|_E=\rho(f^*)\}$ are the corresponding quasi-Banach function space and r.i. quasi-Banach function space, respectively. Now, we are interested in the boundedness of the averaging operators $H$ from Definition 14 from $E$ to $B$.

First, we make the following remark. In the case $d=1$, that is if the averaging operator $H$ is either the operator $H^{(U,W)}$ or $H_{(V,W)}$, then by Lemma 16 the function $Hf$ is non-increasing and hence $Hf=(Hf)^*$. Thus,

$$\|Hf\|_B=\|(Hf)^*\|_B=\|Hf\|_E.$$

In the general case $d\neq 1$, this is not so. This shows the following example. Let $d(t)=2-\sin t$ and $f(v)=1$, $v\in(0,1)$. It is not difficult to show that the function

$$\left(H^{(1,1;d)}f\right)(t)=\frac{1}{t(2-\sin t)}(2t+1-\cos t)$$

is increasing on $(0,1)$ and hence, $\left(H^{(1,1;d)}f\right)^*\neq H^{(1,1;d)}f$.

The next theorem deals with the case $d=1$. It will play a crucial role in our research below. For the case $L=\infty$ it coincides with Theorem 2 from [59]. Theorem 17 from [52] extends Theorem 2 from [59] for the case $L=1$. Remark that [5, Chapter 3, Theorem 5.15, p. 150] is a special case of [59, Theorem 2].

**Theorem 17**

(i) *If $0<W<\infty$, then the operator $H^{(U,W)}$ is bounded from $E$ to $B$ if and only if $p_E>U$.*
(ii) *The operator $H^{(U,\infty)}$ is bounded from $E$ to $B$ if $p_E>U$.*
(iii) *If $0<W<\infty$, then the operator $H_{(V,W)}$ is bounded from $E$ to $B$ if and only if $q_E<V\le\infty$.*
(iv) *The operator $H_{(V,\infty)}$ is bounded from $E$ to $B$ if $q_E<V<\infty$.*

Reverse implications are not true in parts (ii) and (iv). See [59]. We would like to emphasize that if $0<W<\infty$ this theorem gives necessary and sufficient conditions for the operators $H^{(U,W)}$ and $H_{(V,W)}$ to be bounded from a r.i. quasi-Banach function space to corresponding quasi-Banach function space in terms of its Boyd indices.

To consider the general case $d \in SV$, we need the following lemma.

**Lemma 18** *Let* (2) *holds.*

(i) *If the operator $H^{(U,W)}$ is bounded from $E$ to $B$, then the operators $H^{(U_1,W;d)}$ are also bounded from $E$ to $B$ for all $0 < U_1 < U$.*

(ii) *If $0 < V < \infty$ and the operator $H_{(V,W)}$ is bounded from $E$ to $B$, then the operators $H_{(V_1,W;d)}$ are also bounded from $E$ to $B$ for all $V < V_1 \leq \infty$.*

*Proof* We will prove only the first statement. The second one can be proven in a similar way. For $\varepsilon = \frac{1}{U_1} - \frac{1}{U}$, we have $\varepsilon > 0$ and the function $s^{\varepsilon} d(s)$ is equivalent to an increasing function. Hence,

$$\left(H^{(U_1,W;d)}f\right)(t) = \left(t^{\frac{1}{U_1}} d(t)\right)^{-1} \left\| s^{\frac{1}{U}-\frac{1}{W}} s^{\varepsilon} d(s) f^*(s) \right\|_{W,(0,t)}.$$

$$\prec \left(t^{\frac{1}{U_1}} d(t)\right)^{-1} t^{\varepsilon} d(t) \left\| s^{\frac{1}{U}-\frac{1}{W}} f^*(s) \right\|_{W,(0,t)}$$

$$= t^{-\frac{1}{U}} \left\| s^{\frac{1}{U}-\frac{1}{W}} f^*(s) \right\|_{W,(0,t)} = \left(H^{(U,W)}f\right)(t).$$

By the lattice property (P2) and because the operator $H^{(U,W)}$ is bounded from $E$ to $B$, we conclude that

$$\left\| H^{(U_1,W;d)} f \right\|_B \prec \left\| H^{(U,W)} f \right\|_B \prec \|f\|_E.$$

Thus, the operator $H^{(U_1,W;d)}$ is also bounded from $E$ to $B$. ■

**Lemma 18** and **Theorem 17** allow us to prove the main result of this subsection.

**Theorem 19** (Cf. [57, Theorem 1, p 2].) *Let* (2) *holds.*

(i) *If $p_E > U$, then the operator $H^{(U,W;d)}$ is bounded from $E$ to $B$.*

(ii) *If $q_E < V \leq \infty$, then the operator $H_{(V,W;d)}$ is bounded from $E$ to $B$.*

*Proof* We prove the assertion (i). The assertion (ii) can be proved in a similar way. Let $U_1$ be such that $U < U_1 < p_E$. Theorem 17 (i) and (ii) implies that the operator $H^{(U_1,W)}$ is bounded from $E$ to $B$. Because $U < U_1$, by Lemma 18 (i), we conclude that the operator $H^{(U,W;d)}$ is also bounded from $E$ to $B$. ■

*Remark 20* Let $-\infty < V < 0$ and $0 < W \leq \infty$. Formally, the operators

$$\left(H_{(V,W;d)}f\right)(t) = \left(t^{\frac{1}{V}} d(t)\right)^{-1} \left\| s^{\frac{1}{V}-\frac{1}{W}} d(s) f^*(s) \right\|_{W,(t,L)}$$

can also be considered for this case. Because $d \in SV$, $f^*$ is non-increasing and $\frac{1}{V} < 0$ and using Lemma 2 (iii), we see that

$$\left(H_{(V,W;d)}f\right)(t) \leq \left(t^{\frac{1}{V}} d(t)\right)^{-1} f^*(t) \left\| s^{\frac{1}{V}-\frac{1}{W}} d(s) \right\|_{W,(t,\infty)}$$

$$\sim \left(t^{\frac{1}{V}} d(t)\right)^{-1} f^*(t) t^{\frac{1}{V}} d(t) = f^*(t).$$

On the contrary, according to Lemma 2 (v) we have

$$\left(H_{(V,W;d)}f\right)(t) \geq \left(t^{\frac{1}{V}} d(t)\right)^{-1} \left\| s^{\frac{1}{V}-\frac{1}{W}} d(s) f^*(s) \right\|_{W,(t,2t)}$$

$$\geq \left(t^{\frac{1}{V}}d(t)\right)^{-1} f^*(2t) \left\| s^{\frac{1}{V}-\frac{1}{W}}d(s) \right\|_{W,(t,2t)}$$
$$\sim \left(t^{\frac{1}{V}}d(t)\right)^{-1} f^*(2t)t^{\frac{1}{V}}d(t) = f^*(2t).$$

Thus,
$$f^*(2t) \prec H_{(V,W;d)}f(t) \prec f^*(t).$$
Hence, by the lattice property (P2) and Lemma 12, we observe that $\forall f \in E$
$$\|f\|_E = \|f^*(t)\|_E \sim \|f^*(2t)\|_E = \|f^*(2t)\|_B \prec \left\|H_{(V,W;d)}f\right\|_B$$
and
$$\left\|H_{(V,W;d)}f\right\|_B \prec \|f^*\|_B = \|f\|_E.$$
So,
$$\left\|H_{(V,W;d)}f\right\|_B \sim \|f\|_E.$$
Therefore, the expression $\left\|H_{(V,W;d)}f\right\|_B$ $(-\infty < V < 0)$ gives an equivalent quasi-norm on each r.i. quasi-Banach function space $E$ on $(0, L)$. We will not further include this trivial case in our results about equivalent quasi-norms on r.i. quasi-Banach function spaces.

## 3 Main results

Next in this section, we also assume that $\rho$ is a function quasi-norm on $(0, L)$ $(L \in \{1, \infty\})$ and $\{B, \|f\|_B = \rho(|f|)\}$ and $\{E, \|f\|_E = \rho(f^*)\}$ are the corresponding quasi-Banach function space and r.i. quasi-Banach function space, respectively. $0 < p_E \leq q_E \leq \infty$ are the Boyd indices of the space $E$. Additionally, we assume that (2) holds. Lemma 15 and Theorem 19 allows us to prove our main result concerning equivalent quasi-norms on r.i. quasi-Banach function spaces.

**Theorem 21** *If* $0 < U < p_E$, *then for all* $f \in \mathfrak{M}(0, L)$
$$\|f\|_E \sim \left\|H^{(U,W;d)}f\right\|_B.$$
*Moreover,* $\left\|H^{(U,W;d)}f\right\|_B$ *is an equivalent quasi-norm on* $E$*, and* $E$ *equipped with this quasi-norm is also a r.i. quasi-Banach function space.*

(i) *If* $q_E < V \leq \infty$, *then for all* $f \in \mathfrak{M}(0, L)$
$$\|f\|_E \sim \left\|H_{(V,W;d)}f\right\|_B.$$
*Moreover,* $\left\|H_{(V,W;d)}f\right\|_B$ *is an equivalent quasi-norm on* $E$*, and* $E$ *equipped with this quasi-norm is also a r.i. quasi-Banach function space.*

*Proof* First, suppose that $H$ is one of the averaging operators $H^{(U,W;d)}$, then by Lemma 15 (vii), for all $f \in \mathfrak{M}(0, L)$ and $t > 0$
$$f^*(t) \prec (Hf)(t).$$
The lattice property (P2) implies now
$$\|f\|_E = \|f^*\|_B \prec \|Hf\|_B. \tag{3}$$
If $H$ is one of the averaging operators $H_{(V,W;d)}$, then by Lemma 15 (viii), we observe similarly that
$$\|f\|_E = \|D_{2^{-1}}D_2 f\|_E \leq \|D_{2^{-1}}\| \|f^*(2t)\|_E \sim \|f^*(2t)\|_B \prec \|Hf\|_B$$
and we also get the estimate (3). Thus, we only need to prove the inverse inequalities and the rest of the statements under the assumption $\|f\|_E < \infty$.

We continue the proof of (i). The part (ii) can be considered analogously. Let $H$ is one of the averaging operators $H^{(U,W;d)}$. Because $0 < U < p_E$, by Theorem 19 (i), we have

$$\|Hf\|_B \prec \|f\|_E.$$

Thus, considering (3),

$$\|f\|_E \sim \|Hf\|_B.$$

This means that for some positive constants $c$ and $C$ for all $f \in \mathfrak{M}(0,L)$, it holds

$$c\|f\|_E \le \|Hf\|_B \le C\|f\|_E. \quad (4)$$

As the next step, we show that $\|Hf\|_B$ is a r.i. quasi-norm on $E$. (See Definition 5.) To prove (P1), it suffices to prove the quasi-triangle inequality. By **Lemma 15** (vi), using Lemma 12 and (4) we have

$$\begin{aligned}
\|(H(f+g))(t)\|_B &\le \left\| C_d \max\left\{1, 2^{\frac{1-W}{W}}\right\}\left\{(Hf)\left(\frac{t}{2}\right) + (Hg)\left(\frac{t}{2}\right)\right\}\right\|_B \\
&\le C_B C_d \max\left\{1, 2^{\frac{1-W}{W}}\right\}\left\{\left\|(Hf)\left(\frac{t}{2}\right)\right\|_B + \left\|(Hg)\left(\frac{t}{2}\right)\right\|_B\right\} \\
&\le C C_B C_d \max\left\{1, 2^{\frac{1-W}{W}}\right\}\left\{\left\|f\left(\frac{t}{2}\right)\right\|_E + \left\|g\left(\frac{t}{2}\right)\right\|_E\right\} \\
&\sim C C_B C_d \max\left\{1, 2^{\frac{1-W}{W}}\right\}\{\|f(t)\|_E + \|g(t)\|_E\} \\
&\le \frac{1}{c} C C_B C_d \max\left\{1, 2^{\frac{1-W}{W}}\right\}\{\|(Hf)(t)\|_B + \|(Hg)(t)\|_B\}
\end{aligned}$$

where $C_B$ is the quasi-triangle constant of the space $\{B, \|f\|_B\}$. So, (P1) is met. The lattice property (P2) follows from Lemma 15 (iv) and the lattice property of the quasi-norm $\|f\|_B$. Let $\omega \subset (0,L)$ be a measurable set with finite measure. By (4) and the property (P4) of the space $E$ (see Remark 8), we have

$$\|H\chi_\omega\|_B \le C\|\chi_\omega\|_E < \infty.$$

So, (P4) is met for $\|Hf\|_B$. Finally, we only need to prove the Fatou property (P3). Let $0 \le f_n \nearrow f \ (n = 1,2,\dots)$ and $\|Hf_n\|_B < \infty$. By [5, Chapter 2, Proposition 1.7]:

$$0 \le f_n \nearrow f \qquad \text{implies} \qquad {f_n}^* \nearrow f^*.$$

Let us consider the case $W < \infty$ and the operators $H^{(U,W;d)}$. According to Lévy's theorem on monotonic convergence, for any $0 < t < L$

$$\int_0^t \left[v^{\frac{1}{U}} d(v) {f_n}^*(v)\right]^W \frac{dv}{v} \nearrow \int_0^t \left[v^{\frac{1}{U}} d(v) f^*(v)\right]^W \frac{dv}{v}.$$

Obviously, if $\int_0^{t_0} \left[v^{\frac{1}{U}} d(v) f^*(v)\right]^W \frac{dv}{v} = \infty$ for some $t_0 < L$, then $\int_0^t \left[v^{\frac{1}{U}} d(v) f^*(v)\right]^W \frac{dv}{v} = \infty$ for all $t > t_0$ and $\|Hf\|_B = \infty$. Otherwise,

$$H^{(U,W;d)} f_n \nearrow H^{(U,W;d)} f.$$

This means that for each $0 < t < L$ either $\left(H^{(U,W;d)} f\right)(t) < \infty$ and $\left(H^{(U,W;d)} f_n\right)(t) \nearrow \left(H^{(U,W;d)} f\right)(t)$ or $\left(H^{(U,W;d)} f\right)(t) = \infty$ and $\left(H^{(U,W;d)} f_n\right)(t) \nearrow \infty$. Because the property (P3) holds for the space $E$, by (4), we conclude that

$$\left\|H^{(U,W;d)} f_n\right\|_B \nearrow \left\|H^{(U,W;d)} f\right\|_B.$$

Thus, (P3) holds for the quasi-norm $\left\|H^{(U,W;d)} f\right\|_B$. The case $W = \infty$ can be considered similarly. By definition, $Hf^* = Hf$. Hence, $\{E, \|Hf\|_B\}$, is rearrangement invariant. ■

*Remark 22* Theorem 21 can be used iteratively. For example, let $\{E, \rho(f^*)\}$ be a r.i. quasi-Banach function space such that $p_E > 1$. Then (using $H^{(1,1)} f = f^{**}$) by Theorem 21 (i), we get

$$\|f\|_E \sim \|f^{**}\|_E.$$

If we consider the space $E$ with the equivalent quasi-norm $\|f^{**}\|_E$ we can use Theorem 21 again. So, we get, for example, the following statements.

(i) If $p_E > 1$ and $0 < U < p_E$, then for all $f \in \mathfrak{M}$

$$\|f^{**}\|_E \sim \left\|H^{(U,W;d)} f^{**}\right\|_B.$$

(ii) If $1 < p_E \le q_E < V \le \infty$, then for all $f \in \mathfrak{M}$

$$\|f^{**}\|_E \sim \left\|H_{(V,W;d)} f^{**}\right\|_B.$$

As an example, see Theorem 36 below.

**Theorem 17** and **Theorem 21** allows us to formulate the following criterion.

**Theorem 23** (Cf. [21, Theorem 11].) *Let* $0 < W < \infty$.

(i) *If* $0 < U < \infty$, *then for all* $f \in \mathfrak{M}(0, L)$

$$\|f\|_E \sim \left\|H^{(U,W)} f\right\|_B, \qquad \text{if and only if } p_E > U.$$

*In particular,*

$$\|f\|_E \sim \left\|f^{**}_{(U)}\right\|_B, \qquad \text{if and only if } p_E > U,$$

*and*

$$\|f\|_E \sim \|f^{**}\|_B, \qquad \text{if and only if } p_E > 1.$$

(ii) *If* $0 < V \le \infty$, *then for all* $f \in \mathfrak{M}(0, L)$

$$\|f\|_E \sim \left\|H_{(V,W)} f\right\|_B, \qquad \text{if and only if } q_E < V.$$

*In particular,*

$$\|f\|_E \sim \left\|t^{-1}\int_t^\infty f^*(s)ds\right\|_B, \qquad \text{if and only if } q_E < 1.$$

*Proof* First, recall that by Lemma 16 in all cases under consideration $\|Hf\|_E = \|Hf\|_B$. We prove (i). The statements (ii) can be proved similarly. Let $0 < U < p_E$ and $0 < W < \infty$. Then from Theorem 21 (i), we know that

$$\|f\|_E \sim \left\|H^{(U,W)} f\right\|_B.$$

It remains to prove that if $0 < U, W < \infty$ and $\|f\|_E \sim \left\|H^{(U,W)} f\right\|_B$ for all $f \in \mathfrak{M}(0, L)$, then $U < p_E$. Indeed, $\|f\|_E \sim \left\|H^{(U,W)} f\right\|_B$ implies that the operator $H^{(U,W)}$ is bounded from $E$ to $B$. Thus, by Theorem 17 (i), we conclude that $U < p_E$. ■

## 4 Fundamentals the interpolation theory of operators

The general results from Section 3 are formulated in terms of Boyd indices. They can often be reformulated in terms of parameters of concrete functional scales. This will be done in Section 5. To do this, we will need some facts from the interpolation theory of operators and spaces. In this section we recall some basic constructions and definitions related to this theory. For more details we refer to [5, 6, 69]. In this paper we are mainly concerned with interpolation within the class of quasi-Banach function spaces, while the bounded quasilinear operators are considered as the corresponding morphisms. Hence, all definitions we give in the form we need.

**Definition 24** ([5, Chapter 3, Definition 5.3, p. 143]) Let $T$ be an operator whose domain is some linear subspace of $\mathfrak{M}$ and whose range is contained in $\mathfrak{M}$. Then $T$ is said to be *C-sublinear* or *quasilinear* if there is a constant $C \ge 1$ such that the relations

$$|T(f+g)(t)| \le C(|Tf(t)| + |Tg(t)|) \quad \text{and} \quad |T(\lambda f)(t)| = |\lambda||T(f)(t)|$$

hold a.e for all $f$ and $g$ in the domain of $T$ and for all scalars $\lambda$. If $C = 1$, then $T$ is said to be *sublinear*.

### 4.1 Basic definitions

In the following, $\bar{E} \equiv (E_0, E_1)$ will be a *compatible couple* of quasi-Banach function spaces on $(0, L)$. This means that $E_0$ and $E_1$ are continuously embedded in some common Hausdorff topological vector space. As usual, the intersection and the sum are defined as

$$E_0 \cap E_1 = \{f \in \mathfrak{M}: f \in E_0, f_1 \in E_1 \}$$

and

$$E_0 + E_1 = \{f \in \mathfrak{M}: \exists f_0 \in E_0, f_1 \in E_1: f = f_0 + f_1 \},$$

equipped with the quasi-norms

$$\|f\|_{E_0 \cap E_1} = \max\{\|f\|_{E_0}, \|f\|_{E_1}\}$$

and

$$\|f\|_{E_0 + E_1} = \inf\{\|f_0\|_{E_0} + \|f_1\|_{E_1}; f_0 \in E_0, f_1 \in E_1: f = f_0 + f_1\},$$

respectively. Obviously, $E_0 \cap E_1$ and $E_0 + E_1$ are also (quasi-) Banach function spaces.

Let $\bar{E} = (E_0, E_1)$ and $\bar{F} = (F_0, F_1)$ be two compatible couples of (quasi-) Banach function spaces. A quasi-Banach function space $E$ is called an *intermediate space* for the couple $\bar{E}$ (or *intermediate space between* $E_0$ and $E_1$) if

$$E_0 \cap E_1 \hookrightarrow E \hookrightarrow E_0 + E_1.$$

The space $E$ is called an *interpolation space* between $E_0$ and $E_1$ (or *with respect to* $\bar{E}$) if in addition

$$T: E_0 + E_1 \to E_0 + E_1 \quad \text{implies} \quad T: E \to E$$

for any quasilinear operator $T$.

**Lemma 25** *Let $E_0$, $E_1$ and $E$ be r.i. quasi-Banach function spaces on $(0, L)$. If the r.i. quasi-Banach function space $E$ is an interpolation space between $E_0$ and $E_1$, then for the Boyd indices of these spaces, it holds*

$$p_E \geq \min\{p_{E_0}, p_{E_1}\} \qquad \textit{and} \qquad q_E \leq \max\{q_{E_0}, q_{E_1}\}.$$

*Proof* Let us choose some $0 < U < \min\{p_{E_0}, p_{E_1}\}$ and consider the averaging operator $H^{(U,U)}$. Theorem 17 (i) (part "if") implies

$$H^{(U,U)}: E_i \to E_i \qquad (i = 0,1).$$

Since $E$ is an interpolation space between $E_0$ and $E_1$, we have

$$H^{(U,U)}: E \to E.$$

Thus, by Theorem 17 (i) (part "only if"), we get $U < p_E$. As $U$ is arbitrary in the open interval $(0, \min\{p_{E_0}, p_{E_1}\})$, we conclude that

$$p_E \geq \min\{p_{E_0}, p_{E_1}\}.$$

The estimate for $q_E$ can be proved in a similar way using Theorem 17 (ii). ■

The previous lemma implies

**Corollary 26** *If a r.i. quasi-Banach function space $E$ on $(0, L)$ is an interpolation space between $L_{p_0}$ and $L_{p_1}$ where $0 < p_0 < p_1 \leq \infty$. Then*

$$p_0 \leq p_E \leq q_E \leq p_1.$$

### 4.2 *K*-Method of real interpolation

We will use one of the most important ways of constructing interpolation spaces based on the use of the *Peetre's K-functional*. For an arbitrary quasi-Banach couple $\bar{X} = (X_0, X_1)$ and for each $t > 0$, the Peetre's $K$-functional is defined on $X_0 + X_1$ and given by

$$K(t, x; \bar{X}) \equiv K(t, x) := \inf(\|x_0\|_{X_0} + t\|x_1\|_{X_1}: x = x_0 + x_1, x_i \in X_i).$$

For a fixed $x \in X_0 + X_1$ the function $t \rightsquigarrow K(t,x)$ is continuous, non-decreasing, concave and non-negative on $(0,\infty)$, the function $t \rightsquigarrow t^{-1}K(t,x)$ is continuous and non-increasing. Additional properties of the $K$-functional can be found, for example, in [5, 6, 69]. Below we assume that they are known to the reader. Obviously, we have

$$K(t,x) \le K(\lambda t,x) \le \lambda K(t,x) \quad (\lambda \ge 1). \tag{5}$$

**Definition 27** ([44]) Let $0 \le \theta \le 1$, $0 < r \le \infty$ and $b \in SV$. We put

$$\bar{X}_{\theta,r;b} \equiv (X_0,X_1)_{\theta,r;b} := \left\{x \in X_0 + X_1 : \|x\|_{\theta,r;b} = \left\|u^{-\theta-1/r}b(u)K(u,x)\right\|_{r,(0,\infty)} < \infty\right\}.$$

If $b = 1$ this definition recovers the classical (Lions–Peetre) interpolation space $\bar{X}_{\theta,r}$.

**Lemma 28** ([44, Proposition 2.5]) *$\bar{X}_{\theta,r;b}$ is a (quasi-) Banach space. $X_0 \cap X_1 \hookrightarrow \bar{X}_{\theta,r;b} \hookrightarrow X_0 + X_1$ if and only if one of the following conditions is satisfied*:

(i) $0 < \theta < 1$,

(ii) $\theta = 0$ *and* $\left\|u^{-1/r}b(u)\right\|_{r,(1,\infty)} < \infty$,

(iii) $\theta = 1$ *and* $\left\|u^{-1/r}b(u)\right\|_{r,(0,1)} < \infty$.

*In these cases, $\bar{X}_{\theta,r;b}$ is an interpolation space between $X_0$ and $X_1$. If none of these conditions hold, then $\bar{X}_{\theta,r;b} = \{0\}$.*

# 5 Applications

## 5.1 Lorentz–Karamata spaces

We begin with the Lorentz–Karamata spaces. They were introduced in 2000 by Edmunds, Kerman and Pick in [27] and form an important scale of spaces. It contains, e.g., Lebesgue spaces $L_p$, Lorentz space $L_{p,r}$, Lorentz–Zygmund and the generalized Lorentz–Zygmund spaces. For further information about Lorentz–Karamata spaces we refer to, e.g., [44, 65]. They have found many different important applications in analysis, see [1, 4, 5, 20, 26, 30, 44, 64, 65] and the references therein.

In this subsection, we generally assume that $0 < p,r \le \infty$, $b \in SV(0,L)$ and (2) holds. We introduce the mapping $\rho_{p,r;b}: \mathfrak{M}^+(0,L) \to [0,\infty]$ by the formula (Cf. Example 10.)

$$\rho_{p,r;b}(f) = \left\|t^{\frac{1}{p}-\frac{1}{r}}b(t)f(t)\right\|_{r,(0,L)}.$$

Obviously, it is a function quasi-norm. The corresponding quasi-Banach function space we denote by $B_{p,r;b}$:

$$B_{p,r;b} \equiv B_{p,r;b}(0,L) = \left\{f \in \mathfrak{M}(0,L): \ \|f\|_{B_{p,r;b}} := \rho_{p,r;b}(|f|) < \infty\right\}.$$

**Definition 29** ([65, Definition 3.1]) The Lorentz–Karamata spaces $L_{p,r;b} \equiv L_{p,r;b}(0,L)$ and $L_{(p,r;b)} \equiv L_{(p,r;b)}(0,L)$ are the sets of all $f \in \mathfrak{M}(0,L)$ such that

$$\|f\|_{p,r;b} := \|f^*\|_{B_{p,r;b}} = \left\|t^{\frac{1}{p}-\frac{1}{r}}b(t)f^*(t)\right\|_{r,(0,L)} < \infty$$

and

$$\|f\|_{(p,r;b)} := \|f^{**}\|_{B_{p,r;b}} = \left\|t^{\frac{1}{p}-\frac{1}{r}}b(t)f^{**}(t)\right\|_{r,(0,L)} < \infty,$$

respectively.

If $b = 1$, we skip the character $b$ in the above notations. It ttis case we get the Lorentz spaces. The following two lemmas establish the basic properties of the Lorentz–Karamata spaces $L_{p,r;b}$ that we need.

**Lemma 30** ([65, Proposition 3.2 and Proposition 3.4])

(i) *The space* $L_{p,r;b}(0,L)$ *is a r.i. quasi-Banach function space, if and only if one of the following conditions is satisfied*:

$$\begin{cases} 0 < p < \infty; \\ p = \infty \text{ and } \left\| t^{-\frac{1}{r}} b(t) \right\|_{r,(0,1)} < \infty. \end{cases}$$

*Furthermore, if those conditions are not satisfied, then* $L_{p,r;b}(0,L) = \{0\}$.

(ii) *The space* $L_{(p,r;b)}(0,\infty)$ *is a r.i. quasi-Banach function space, if and only if one of the following conditions is satisfied*:

$$\begin{cases} 1 < p < \infty; \\ p = 1 \text{ and } \left\| t^{-\frac{1}{r}} b(t) \right\|_{r,(1,\infty)} < \infty; \\ p = \infty \text{ and } \left\| t^{-\frac{1}{r}} b(t) \right\|_{r,(0,1)} < \infty. \end{cases} \tag{6}$$

*Furthermore, if none of these conditions is satisfied, then* $L_{(p,r;b)}(0,\infty) = \{0\}$.

Let us consider the spaces $L_{(p,r;b)}(0,1)$. The following lemma can be proved using standard methods. See Lemma 3.15 in [64].

**Lemma 31** (Cf. [64, Lemma 3.15].) *Let one of the following conditions be satisfied*:

$$\begin{cases} 0 < p < 1; \\ p = 1 \; and \; \left\| t^{-\frac{1}{r}} b(t) \right\|_{r,(0,1)} < \infty. \end{cases}$$

*Then*

$$L_{(p,r;b)}(0,1) = L_1.$$

Taking (6) into account, we will consider the spaces $L_{(p,r;b)}(0,L)$ under following assumptions:

$$\begin{cases} 1 < p < \infty; \\ p = 1 \text{ and } \left\| t^{-\frac{1}{r}} b(t) \right\|_{r,(0,1)} = \infty, if \; L = 1; \\ p = 1 \text{ and } \left\| t^{-\frac{1}{r}} b(t) \right\|_{r,(1,\infty)} < \infty, if \; L = \infty; \\ p = \infty \text{ and } \left\| t^{-\frac{1}{r}} b(t) \right\|_{r,(0,1)} < \infty. \end{cases} \tag{7}$$

We will also need some interpolation formulae of Lorentz–Karamata spaces.

**Lemma 32**

(i) *If* $0 < p_0 < p < p_1 \leq \infty$, *then there exist* $0 < \theta < 1$ *and* $\hat{b} \in SV$ *such that*

$$L_{p,r;b} = \left(L_{p_0}, L_{p_1}\right)_{\theta,r;\hat{b}}.$$

(ii) *Let* $\left\| t^{-\frac{1}{r}} b(t) \right\|_{r,(0,1)} < \infty$. *Then for each* $0 < p_0 < \infty$ *there exists* $\hat{b} \in SV$ *such that*

$$L_{\infty,r;b} = \left(L_{p_0}, L_\infty\right)_{1,r;\hat{b}}.$$

(iii) *Let* $\left\| u^{-1/r} b(u) \right\|_{r,(0,1)} = \infty$ *if* $L = 1$ *and* $\left\| u^{-1/r} b(u) \right\|_{r,(1,\infty)} < \infty$ *if* $L = \infty$. *Then for each* $1 < p_0 < \infty$ *there exists* $\hat{b} \in SV$ *such that*

$$L_{(1,r;b)} = \left(L_1, L_{p_0}\right)_{0,r;\hat{b}}.$$

*Proof* The statements (i) and (ii) follow from [44, Corollary 5.3, Lemma 5.5] and [23, Theorem 5.1], respectively. The statement (iii) follows from [44, (6.7), case $\theta = 0$] and [44, Theorem 3.5, (3.20)]. ∎

Now, Lemma 25, Lemma 28 and Lemma 32 allow us to calculate the Boyd indices for the spaces $L_{p,r;b}$. Note that this result is proved as a corresponding part of Theorem 3.15 in [65] for the case $L = \infty$. Cf. also [48, Example 1]. Here we give a different proof.

**Corollary 33** *Let $p, r, b$ be such that $L_{p,r;b}$ is non-trivial.* (See Lemma 30 (i).) *Then*

$$p_{L_{p,r;b}} = q_{L_{p,r;b}} = p.$$

*Proof* Case $p < \infty$. By Lemma 32 (i), for any $p_0$ and $p_1$ such that $0 < p_0 < p < p_1 \le \infty$,

$$L_{p,r;b} = \left(L_{p_0}, L_{p_1}\right)_{\theta,r;\hat{b}}$$

for some $0 < \theta < 1$ and $\hat{b} \in SV$. Lemma 25 and Lemma 28 (i) immediately imply that

$$p_{L_{p,r;b}} \ge \sup_{0<p_0<p} p_0 = p$$

and

$$q_{L_{p,r;b}} \le \inf_{p<p_1} p_1 = p.$$

Case $p = \infty$. The space $L_{\infty,r;b}$ is non-trivial if and only if $\left\| t^{-\frac{1}{r}} b(t) \right\|_{r,(0,1)} < \infty$. By Lemma 32 (ii), for any $p_0$ such that $0 < p_0 < \infty$,

$$L_{\infty,r;b} = \left(L_{p_0}, L_\infty\right)_{1,r;\hat{b}}$$

for some $\hat{b} \in SV$. Thus, Lemma 25 and Lemma 28 (iii) immediately imply

$$p_{L_{\infty,r;b}} \ge \sup_{0<p_0<\infty} p_0 = \infty.$$

Hence, $q_{L_{\infty,r;b}} = \infty$. ∎

Now, we can formulate results on equivalent quasi-norms on Lorentz–Karamata spaces. Combining **Corollary 33** and **Theorem 21**, we get the following theorem.

**Theorem 34** *Let $p, r, b$ be such that $L_{p,r;b}$ is non-trivial.*

(i) (Cf. [14, Theorem 5.2].) *If $0 < U < p \le \infty$, then for all $f \in \mathfrak{M}(0, L)$*

$$\|f\|_{p,r;b} \sim \left\| H^{(U,W;d)} f \right\|_{B_{p,r;b}}.$$

*Moreover, $\left\| H^{(U,W;d)} f \right\|_{B_{p,r;b}}$ is an equivalent r.i. quasi-norm on $L_{p,r;b}$.*

(ii) *If $p < V \le \infty$,* then for all $f \in \mathfrak{M}(0, L)$

$$\|f\|_{p,r;b} \sim \left\| H_{(V,W;d)} f \right\|_{B_{p,r;b}}.$$

*Moreover, $\left\| H_{(V,W;d)} f \right\|_{B_{p,r;b}}$ is an equivalent r.i. quasi-norm on $L_{p,r;b}$.*

*Remark 35*

(i) **Theorem 34** (i) is equivalent to [29, Theorem 3.1].

(ii) Theorem 34 (ii), in conjunction with Remark 20, improves [29, Theorem 3.2].

Using additionally *Remark 22*, we arrive at follows

**Theorem 36** *Let $p, r, b$ be such that* (7) *holds.*

(i) *If $1 < p \le \infty$ and $0 < U < p$, then for all $f \in \mathfrak{M}(0, L)$*

$$\|f\|_{(p,r;b)} \sim \left\| H^{(U,W;d)} f^{**} \right\|_{B_{p,r;b}}.$$

*Moreover, $\left\| H^{(U,W;d)} f^{**} \right\|_{B_{p,r;b}}$ is an equivalent r.i. quasi-norm on $L_{(p,r;b)}$.*

(ii) *If* $1 < p < V \leq \infty$, *then for all* $f \in \mathfrak{M}(0, L)$

$$\|f\|_{(p,r;b)} \sim \|H_{(V,W;d)} f^{**}\|_{B_{p,r;b}}.$$

*Moreover,* $\|H_{(V,W;d)} f^{**}\|_{B_{p,r;b}}$ *is an equivalent r.i. quasi-norm on* $L_{(p,r;b)}$.

*Remark 37*

(i) **Theorem 36** (i) is equivalent to [29, Theorem 3.5]. See also [28, Theorems 1.2 and 1.3].

(ii) **Theorem 36** (ii), in conjunction with Remark 20, improves [29, Theorem 3.6]. See also [28, Theorems 1.1 and 1.4].

*Remark 38* For those interested in equivalent quasi-norms in Lorentz–Karamata spaces, we recommend the paper [33] by P. Fernández-Martínez and T. Signes. From a mathematical point of view, Theorem 34 and [33, Theorems 3.6 and 3.7] demonstrate a large overlap and symmetric difference. Fernández-Martínez and Signes work with the rich scale of Lorentz–Karamata type spaces $\left\|t^{\frac{1}{p}} b(t) f^*(t)\right\|_{\tilde{E}}$ based on the family of r.i. Banach function spaces $E$. However, the spaces $L_{p,r;b}$ with $0 < r < 1$ do not belong to this scale.

Thanks to Theorem 23, some statements of Theorem 34 can be converted into criteria for the averaging operators $H^{(U,W)}$ and $H_{(V,W)}$.

**Theorem 39** *Let* $p, r, b$ *be such that* $L_{p,r;b}$ *is non-trivial and* $0 < W < \infty$.

(i) $\|f\|_{p,r;b} \sim \|H^{(U,W)} f\|_{B_{p,r;b}}$ *if and only if* $p > U$.

(ii) $\|f\|_{p,r;b} \sim \|H_{(V,W)} f\|_{B_{p,r;b}}$ *if and only if* $p < V$.

In the following corollary, we formulate some particular cases. (Cf. [29, Remark 3.7], [70, Proposition 3], [65, Corollary 3.22] and [24, Theorem 6.2].)

**Corollary 40**

(i) *By Theorem 39* (i) (*with* $W = U$),

$$\|f\|_{p,r;b} \sim \|f^{**}_{(U)}\|_{p,r;b} \quad \text{if and only if } p > U.$$

*In particular*,

$$\|f\|_{p,r;b} \sim \|f\|_{(p,r;b)} \quad \text{if and only if } p > 1.$$

(ii) *By Theorem 39* (ii) (*with* $V = \infty$ *and* $W = 1$),

$$\|f\|_{p,r;b} \sim \left\|t^{\frac{1}{p}-\frac{1}{r}} b(t) \int_t^\infty v^{-1} f^*(v) dv\right\|_{r,(0,\infty)} \quad \text{if and only if } p < \infty.$$

(iii) *By Theorem 39* (ii) (*with* $V = W = 1$)

$$\|f\|_{p,r;b} \sim \left\|t^{\frac{1}{p}-\frac{1}{r}-1} b(t) \int_t^\infty f^*(v) dv\right\|_{r,(0,\infty)} \quad \text{if and only if } p < 1.$$

We can now extend the part of Theorem 3.15 from [65] that concerns $L_{(p,r;b)}$ spaces for the case $L = 1$. (Cf. [48, Example 1].)

**Corollary 41** *Let* $E = L_{(p,r;b)}(0, L)$, *where* $p, r, b$ *fulfil the assumptions* (7). *Then*

$$p_E = q_E = p.$$

*Proof* The case $1 < p \leq \infty$ follows from Corollary 40 (ii) and Corollary 33. The case $p = 1$ follows from Lemma 32 (iii) using logic of the proof of Corollary 33. ■

Using **Lemma 25**, **Corollary 33** and **Corollary 41**, we can extend **Corollary 26**.

**Corollary 42** *Let* $0 < p_0 \le p_1 \le \infty$, $E_0 = L_{p_0,r_0,b_0}(0,L)$ *and* $E_1 = L_{p_1,r_1,b_1}(0,L)$. *Assume that the spaces* $E_0$ *and* $E_1$ *are non-trivial. If a r.i. quasi-Banach function space* $E$ *is an interpolation space between* $E_0$ *and* $E_1$*, then for the Boyd indices* $p_E$ *and* $q_E$ *it holds*

$$p_0 \le p_E \le q_E \le p_1.$$

Note that in [21, Theorem 15], it has been proven that, under certain additional conditions, the converse is also true.

**Corollary 43** *Let* $1 \le p_0 \le p_1 \le \infty$, $E_0 = L_{(p_0,r_0,b_0)}(0,L)$ *and* $E_1 = L_{(p_1,r_1,b_1)}(0,L)$. *Assume that the triples* $(p_0, r_0, b_0)$ *and* $(p_1, r_1, b_1)$ *fulfil the assumptions* (7). *If a r.i. quasi-Banach function space* $E$ *is an interpolation space between* $E_0$ *and* $E_1$*, then for the Boyd indices* $p_E$ *and* $q_E$ *it holds*

$$p_0 \le p_E \le q_E \le p_1.$$

### 5.2 $\mathcal{L}$- and $\mathcal{R}$-Lorentz–Karamata spaces

In this subsection, we generally assume that $0 < p, r, q \le \infty$, $a, b \in SV(0,L)$ and (2) holds. We introduce two mappings $\rho^{\mathcal{L}}_{p,r,b,q,a}$ and $\rho^{\mathcal{R}}_{p,r,b,q,a}$ from $\mathfrak{M}^+(0,L)$ into $[0,\infty]$ by the formulae

$$\rho^{\mathcal{L}}_{p,r,b,q,a}(f) := \left\| t^{-\frac{1}{r}} b(t) \left\| u^{\frac{1}{p}-\frac{1}{q}} a(u) f(u) \right\|_{q,(0,t)} \right\|_{r,(0,L)}$$

and

$$\rho^{\mathcal{R}}_{p,r,b,q,a}(f) := \left\| t^{-\frac{1}{r}} b(t) \left\| u^{\frac{1}{p}-\frac{1}{q}} a(u) f(u) \right\|_{q,(t,L)} \right\|_{r,(0,L)},$$

respectively. Obviously, they are function quasi-norms. The corresponding quasi-Banach function spaces we denote by $B^{\mathcal{L}}_{p,r,b,q,a}$ and $B^{\mathcal{R}}_{p,r,b,q,a}$, respectively:

$$B^{\mathcal{L}}_{p,r,b,q,a} \equiv B^{\mathcal{L}}_{p,r,b,q,a}(0,L) = \left\{ f \in \mathfrak{M}(0,L) \colon \ \|f\|_{B^{\mathcal{L}}_{p,r,b,q,a}} := \rho^{\mathcal{L}}_{p,r,b,q,a}(|f|) < \infty \right\},$$

$$B^{\mathcal{R}}_{p,r,b,q,a} \equiv B^{\mathcal{R}}_{p,r,b,q,a}(0,L) = \left\{ f \in \mathfrak{M}(0,L) \colon \ \|f\|_{B^{\mathcal{R}}_{p,r,b,q,a}} := \rho^{\mathcal{R}}_{p,r,b,q,a}(|f|) < \infty \right\}.$$

Note that these spaces can be considered as weighted Cesàro, Copson or Tandori function spaces. See [50, p. 1162], [43, Definition 12] and [42].

**Definition 44** (See [18, Definition 34] and [24, Definition 50]. Cf. [23, Definition 3.8], [44, (5.21), (5.33)], [15, Definition 7.1], [36].) The spaces $L^{\mathcal{L}}_{p,r,b,q,a}$, $L^{\mathcal{R}}_{p,r,b,q,a}$, $(L)^{\mathcal{L}}_{p,r,b,q,a}$ and $(L)^{\mathcal{R}}_{p,r,b,q,a}$ are the sets of all $f \in \mathfrak{M}(0,L)$ such that

$$\|f\|_{L^{\mathcal{L}}_{p,r,b,q,a}} := \|f^*\|_{B^{\mathcal{L}}_{p,r,b,q,a}} < \infty,$$

$$\|f\|_{(L)^{\mathcal{L}}_{p,r,b,q,a}} := \|f^{**}\|_{B^{\mathcal{L}}_{p,r,b,q,a}} < \infty,$$

$$\|f\|_{L^{\mathcal{R}}_{p,r,b,q,a}} := \|f^*\|_{B^{\mathcal{R}}_{p,r,b,q,a}} < \infty,$$

and

$$\|f\|_{(L)^{\mathcal{R}}_{p,r,b,q,a}} := \|f^{**}\|_{B^{\mathcal{R}}_{p,r,b,q,a}} < \infty,$$

respectively.

We refer to these spaces as $\mathcal{L}$- and $\mathcal{R}$-Lorentz–Karamata spaces, respectively. They naturally arise in reiteration formulae for the limiting cases $\theta = 0$ or $\theta = 1$. In literature (see, [1, 15, 16, 17, 18, 19, 22, 24, 25, 26, 30, 32, 34, 35] and [44, Theorem 31]) similar definitions are given alongside with properties of these spaces. Note that the spaces $L^{\mathcal{L}}_{p,r,b,q,a}$ are special cases of generalized gamma space with double weights [38, Definition 1.2]. In [3, Definition 5.4] they are denoted by $\Gamma(q,r,p;b,a)$. About generalized gamma spaces, see also [2, Definition 2.3]. We will consider only $\mathcal{L}$-spaces. The $\mathcal{R}$-spaces can be considered similarly. We leave this to the reader. Obviously,

$$\|f\|_{L^{\mathcal{L}}_{p,r,b,q,a}} = \left\|t^{\frac{1}{p}-\frac{1}{r}}a(t)b(t)\big(H^{(p,q;a)}f\big)(t)\right\|_{r,(0,L)} = \left\|H^{(p,q;a)}f\right\|_{B_{p,r;ab}}.$$

Using standard arguments, the following two lemmas can be proved.

**Lemma 45** *The space* $L^{\mathcal{L}}_{p,r,b,q,a}$ *is non-trivial if and only if one of the following conditions holds*:

(i) $p=\infty$, $\left\|t^{-\frac{1}{r}}b(t)\right\|_{r,(1,\infty)} < \infty$ *and* $\left\|t^{-\frac{1}{r}}b(t)\left\|u^{-\frac{1}{q}}a(u)\right\|_{q,(0,t)}\right\|_{r,(0,1)} < \infty$,

(ii) $p<\infty$ *and* $\left\|t^{-\frac{1}{r}}b(t)\right\|_{r,(1,\infty)} < \infty$,

*otherwise* $L^{\mathcal{L}}_{p,r,b,q,a}=\{0\}$. *If* $L=1$, *all conditions including* $\infty$ *as a limit of integration can be omitted.*

The spaces $(L)^{\mathcal{L}}_{p,r,b,q,a}$ arise in this paper only for $p>1$. See Corollary 51.

**Lemma 46** Let $p>1$. *The space* $(L)^{\mathcal{L}}_{p,r,b,q,a}$ *is non-trivial if one of the following conditions holds*:

(i) $p=\infty$, $\left\|t^{-\frac{1}{r}}b(t)\right\|_{r,(1,\infty)} < \infty$ *and* $\left\|t^{-\frac{1}{r}}b(t)\left\|u^{-\frac{1}{q}}a(u)\right\|_{q,(0,t)}\right\|_{r,(0,1)} < \infty$,

(ii) $1<p<\infty$ *and* $\left\|t^{-\frac{1}{r}}b(t)\right\|_{r,(1,\infty)} < \infty$,

*otherwise* $(L)^{\mathcal{L}}_{p,r,b,q,a}=\{0\}$. *If* $L=1$, *all conditions including* $\infty$ *as a limit of integration can be omitted.*

We will require that the parameters of the spaces $L^{\mathcal{L}}_{p,r,b,q,a}$ and $(L)^{\mathcal{L}}_{p,r,b,q,a}$ under consideration satisfies conditions from Lemma 45 or Lemma 46, respectively. The next lemma characterizes the spaces $L^{\mathcal{L}}_{p,r,b,q,a}$ as appropriate limiting interpolation spaces.

**Lemma 47** *Let the space* $L^{\mathcal{L}}_{p,r,b,q,a}$ *is non-trivial.*

(i) *If* $0<p<p_1<\infty$, *then there exists* $\hat{b}\in SV$ *such that*

$$L^{\mathcal{L}}_{p,r,b,q,a} = \left(L_{p,q;a}, L_{p_1}\right)_{0,r;\hat{b}}.$$

(ii) *If* $0<p_0<\infty$, *then there exists an interpolation functor* $\mathcal{F}$ *such that*

$$L^{\mathcal{L}}_{\infty,r,b,q,a} = \mathcal{F}\left(L_{p_0}, L_\infty\right).$$

*Proof* The assertion (i) follows from [44, Theorem 5.7]. The statement (ii) follows from [24, Theorem 5.2]. ■

The following corollary follows directly from Lemma 25 and **Lemma 47**. The proof is like the proof of **Corollary 33**.

**Corollary 48** *Let the space* $L^{\mathcal{L}}_{p,r,b,q,a}$ *is non-trivial. Put* $E=L^{\mathcal{L}}_{p,r,b,q,a}$. *Then*

$$p_E = q_E = p.$$

**Theorem 21** and **Corollary 48** immediately imply

**Theorem 49** *Let the space* $L^{\mathcal{L}}_{p,r,b,q,a}$ *is non-trivial.*

(i) *If* $0<U<p$, *then for all* $f\in\mathfrak{M}$

$$\|f\|_{L^{\mathcal{L}}_{p,r,b,q,a}} \sim \left\|H^{(U,W;d)}f\right\|_{B^{\mathcal{L}}_{p,r,b,q,a}}.$$

(ii) *If* $p < V \le \infty$, *then for all* $f \in \mathfrak{M}$

$$\|f\|_{L^{\mathcal{L}}_{p,r,b,q,a}} \sim \left\|H_{(V,W;d)}\right\|_{B^{\mathcal{L}}_{p,r,b,q,a}}.$$

Thanks to Theorem 23, for the averaging operators $H^{(U,W)}$ or $H_{(V,W)}$, some statements of Theorem 49 can be converted into criteria.

**Theorem 50** *Let the space* $L^{\mathcal{L}}_{p,r,b,q,a}$ *is non-trivial and* $0 < W < \infty$. *Then*

(i) $\|f\|_{L^{\mathcal{L}}_{p,r,b,q,a}} \sim \left\|H^{(U,W)}f\right\|_{B^{\mathcal{L}}_{p,r,b,q,a}}$ *if and only if* $p > U$.

(ii) $\|f\|_{L^{\mathcal{L}}_{p,r,b,q,a}} \sim \left\|H_{(V,W)}\right\|_{B^{\mathcal{L}}_{p,r,b,q,a}}$ *if and only if* $p < V$.

In the following corollary, we formulate two particular cases.

**Corollary 51** *Let the space* $L^{\mathcal{L}}_{p,r,b,q,a}$ *is non-trivial. By Theorem 50* (i) (*with* $W = U$)

$$\|f\|_{L^{\mathcal{L}}_{p,r,b,q,a}} \sim \left\|f^{**}_{(U)}\right\|_{B^{\mathcal{L}}_{p,r,b,q,a}} \qquad \text{if and only if } p > U.$$

*In particular* (Cf. [23, Theorem 6.3].),

$$\|f\|_{L^{\mathcal{L}}_{p,r,b,q,a}} \sim \|f^{**}\|_{B^{\mathcal{L}}_{p,r,b,q,a}} = \|f\|_{(L)^{\mathcal{L}}_{p,r,b,q,a}} \qquad \text{if and only if } p > 1.$$

### 5.3 Grand and small Lorentz spaces involving slowly varying functions

The classical grand Lebesgue spaces $L^{p)}$ $(1 < p < \infty)$ were introduced in 1992 by T. Iwaniec and C. Sbordone in [46] in connection with the study of integrability properties of Jacobian. The classical small Lebesgue space $L^{(p}$ was presented by A. Fiorenza in [37] as the associate space of $L^{p)}$. Grand and small Lebesgue spaces and their generalizations - grand and small Lorentz spaces $L^{p),*}_{*}$ and $L^{(p,*}_{*}$ (here "$*$" stands for further parameters) - found numerous important applications, for example, in boundary value problems, embeddings of Besov spaces, PDEs, boundedness of singular integral operators and Riemann–Liouville fractional integral operators, maximal operators and, more generally, quasilinear operators, to describe the behavior of the Fourier coefficients of functions that belong to spaces which are "very close" to $L_2$. For further details, we refer to [4, 17, 20, 26, 34, 35, 38, 39, 41] and the references therein.

These spaces are often defined on bounded domains in $\mathbb{R}^n$. Hence, in this subsection we consider functions from $\mathfrak{M}(0,1)$. Additionally, we generally assume that $0 < p, q, r \le \infty$, $b \in SV(0,1)$ and (2) holds. We introduce two mappings $\rho^{\mathcal{L}}_{p,r,b,q}$ and $\rho^{\mathcal{R}}_{p,r,b,q}$ from $\mathfrak{M}^+(0,1)$ into $[0,\infty]$ by the formulae

$$\rho^{\mathcal{L}}_{p,r,b,q}(f) := \left\| t^{-\frac{1}{r}} b(t) \left\| u^{\frac{1}{p}-\frac{1}{q}} f(u) \right\|_{q,(0,t)} \right\|_{r,(0,1)}$$

and

$$\rho^{\mathcal{R}}_{p,r,b,q}(f) := \left\| t^{-\frac{1}{r}} b(t) \left\| u^{\frac{1}{p}-\frac{1}{q}} f(u) \right\|_{q,(t,1)} \right\|_{r,(0,1)},$$

respectively. Obviously, they are function quasi-norms. The corresponding quasi-Banach function spaces we denote by $B^{(p,q,r}_{b}$ and $B^{p),q,r}_{b}$, respectively:

$$B^{(p,q,r}_{b} = \left\{ f \in \mathfrak{M}(0,1) \colon \ \|f\|_{B^{(p,q,r}_{b}} := \rho^{\mathcal{L}}_{p,r,b,q}(|f|) < \infty \right\},$$

$$B^{p),q,r}_{b} = \left\{ f \in \mathfrak{M}(0,1) \colon \ \|f\|_{B^{p),q,r}_{b}} := \rho^{\mathcal{R}}_{p,r,b,q}(|f|) < \infty \right\}.$$

**Definition 52** The small Lorentz spaces $L_b^{(p,q,r}$ and $(L)_b^{(p,q,r}$ and the grand Lorentz spaces $L_b^{p),q,r}$ and $(L)_b^{p),q,r}$ we define as the sets of all $f \in \mathfrak{M}(0,1)$ such that

$$\|f\|_{L_b^{(p,q,r}} := \|f^*\|_{B_b^{(p,q,r}} < \infty,$$

$$\|f\|_{(L)_b^{(p,q,r}} := \|f^{**}\|_{B_b^{(p,q,r}} < \infty,$$

$$\|f\|_{L_b^{p),q,r}} := \|f^*\|_{B_b^{p),q,r}} < \infty,$$

and

$$\|f\|_{(L)_b^{p),q,r}} := \|f^{**}\|_{B_b^{p),q,r}} < \infty,$$

respectively.

The spaces $L_b^{(p,q,r}$ are introduced in [4]. See also [38, 39, 41] and references therein. For the sake of brevity, we will consider only small Lorentz spaces. Grand Lorentz spaces can be considered similarly. We leave this to the reader. Note that it makes sense to consider the small Lorentz space $L_b^{(p,q,r}$ only under assumptions

$$0 < p < \infty, 0 < q, r \le \infty, b \in SV \text{ and } \left\|t^{-\frac{1}{r}} b(t)\right\|_{r,(0,1)} = \infty. \tag{8}$$

Otherwise, it is either trivial or coincides with the space $L_{p,q}(0,1)$.

Because $f^{**}(u) \ge f^*(u)$, always

$$(L)_b^{(p,q,r} \hookrightarrow L_b^{(p,q,r}.$$

It is clear, that $L_b^{(p,q,r} = L_{p,r,b,q,1}^{\mathcal{L}}(0,1)$ and $(L)_b^{(p,q,r} = (L)_{p,r,b,q,1}^{\mathcal{L}}(0,1)$. Hence, Theorem 49, Theorem 50 and Corollary 51 immediately imply the following three statements.

**Theorem 53** *Let conditions* (8) *be satisfied.*

(i) *If* $0 < U < p$, *then for all* $f \in \mathfrak{M}(0,1)$

$$\|f\|_{L_b^{(p,q,r}} \sim \left\|H^{(U,W;d)} f\right\|_{B_b^{(p,q,r}}.$$

(ii) *If* $p < V$, *then for all* $f \in \mathfrak{M}(0,1)$

$$\|f\|_{L_b^{(p,q,r}} \sim \left\|H_{(V,W;d)}\right\|_{B_b^{(p,q,r}}.$$

**Theorem 54** *Let conditions* (8) *be satisfied and* $0 < W < \infty$.

(i) $\|f\|_{L_b^{(p,q,r}} \sim \left\|H^{(U,W)} f\right\|_{B_b^{(p,q,r}}$ *if and only if* $p > U$.

(ii) $\|f\|_{L_b^{(p,q,r}} \sim \left\|H_{(V,W)}\right\|_{B_b^{(p,q,r}}$ *if and only if* $p < V$.

**Corollary 55** *Let conditions* (8) *be satisfied. By Theorem 54* (i) (*with* $W = U$)

$$\|f\|_{L_b^{(p,q,r}} \sim \left\|f_{(U)}^{**}\right\|_{B_b^{(p,q,r}} \quad \text{if and only if } p > U.$$

*In particular* (Cf. [23, Theorem 6.3].),

$$\|f\|_{L_b^{(p,q,r}} \sim \|f\|_{(L)_b^{(p,q,r}} \quad \text{if and only if } p > 1.$$

We will need the following technical lemma which can be easily proved.

**Lemma 56** *For all non-negative and non-increasing functions* $f \in \mathfrak{M}(0,1)$

$$\left\|t^{\frac{1}{p}-\frac{1}{q}} b(t) f(t)\right\|_{q,(0,1)} \sim \left\|t^{\frac{1}{p}-\frac{1}{q}} b(t) f(t)\right\|_{q,\left(0,\frac{1}{2}\right)}.$$

The following lemma shows why it makes sense to refer to the spaces $L_b^{(p,q,r}$ as “small Lorentz spaces.”

**Lemma 57** (Cf. [4, Remark 2.5].) *Let conditions* (8) *be satisfied, then*

$$L_{x,y;a} \hookrightarrow L_b^{(p,q,r} \hookrightarrow L_{p,q}$$

*for all* $x, y, a$ *such that* $p < x \leq \infty$, $0 < y \leq \infty$ *and* $a \in SV$.

*Proof* Step 1. Using Lemma 56, we get for all $f \in L_b^{(p,q,r}$

$$\|f\|_{L_b^{(p,q,r}} \geq \left\| t^{-\frac{1}{r}} b(t) \left\| u^{\frac{1}{p}-\frac{1}{q}} f^*(u) \right\|_{q,(0,t)} \right\|_{r,\left(\frac{1}{2},1\right)}$$

$$\geq \left\| t^{-\frac{1}{r}} b(t) \right\|_{r,\left(\frac{1}{2},1\right)} \left\| u^{\frac{1}{p}-\frac{1}{q}} f^*(u) \right\|_{q,\left(0,\frac{1}{2}\right)} \sim \left\| u^{\frac{1}{p}-\frac{1}{q}} f^*(u) \right\|_{q,\left(0,\frac{1}{2}\right)} \sim \|f\|_{L_{p,q}}.$$

Thus, $L_b^{(p,q,r} \hookrightarrow L_{p,q}$.

Step 2. Now, we show that $L_{x,q;a} \hookrightarrow L_b^{(p,q,r}$ for all $p < x \leq \infty$ and $a \in SV$. We observe that $\frac{x-p}{px} > 0$. Hence, $\left\| t^{\frac{x-p}{px}-\frac{1}{r}} b(t) \left(a(t)\right)^{-1} \right\|_{r,(0,1)} < \infty$ and $\left(a(u)\right)^{-1} u^{\frac{x-p}{px}}$ is almost increasing. Therefore,

$$\|f\|_{L_b^{(p,q,r}} = \left\| t^{-\frac{1}{r}} b(t) \left\| u^{\frac{1}{x}+\frac{x-p}{px}-\frac{1}{q}} \left(a(u)\right)^{-1} a(u) f^*(u) \right\|_{q,(0,t)} \right\|_{r,(0,1)}$$

$$\prec \left\| t^{\frac{x-p}{px}-\frac{1}{r}} b(t) \left(a(t)\right)^{-1} \left\| u^{\frac{1}{x}-\frac{1}{q}} a(u) f^*(u) \right\|_{q,(0,t)} \right\|_{r,(0,1)}$$

$$\leq \left\| t^{\frac{x-p}{px}-\frac{1}{r}} b(t) \left(a(t)\right)^{-1} \right\|_{r,(0,1)} \left\| u^{\frac{1}{x}-\frac{1}{q}} a(u) f^*(u) \right\|_{q,(0,1)} \sim \|f\|_{L_{x,q;a}}.$$

Step 3. Finally, we show that $L_{x,y;a} \hookrightarrow L_b^{(p,q,r}$ for all $p < x \leq \infty$, $0 < y \leq \infty$ and $a \in SV$. Let $\varepsilon$ is positive and sufficiently small. Then by [65, Theorem 3.17]

$$L_{x,y;a} \hookrightarrow L_{x-\varepsilon,q;a}.$$

Using now Step 2, we arrive at

$$L_{x,y;a} \hookrightarrow L_{x-\varepsilon,q;a} \hookrightarrow L_b^{(p,q,r}. \qquad \blacksquare$$

## References

1. Ahmed, I., Edmunds, D. E., Evans, W. D., Karadzhov, G. E.: Reiteration theorems for the $K$-interpolation method in limiting cases. Math. Nachr. **284**(4) 421‐442 (2011). https://doi.org/10.1002/mana.200810108
2. Ahmed, I., Fiorenza, A., Formica, M. R.: Interpolation of generalized gamma spaces in a critical case. J Fourier Anal. Appl. **28**, 54 (2022). https://doi.org/10.1007/s00041-022-09947-1
3. Ahmed, I., Fiorenza, A., Gogatishvili A.: Holmstedt's formula for the $K$-functional: the limit case $\theta_0=\theta_1$, Math. Nachr. **296**(12), 5474-5492 (2023). https://doi.org/10.1002/mana.202200440
4. Ahmed, I., Fiorenza, A., Hafeez, A.: Some Interpolation Formulae for Grand and Small Lorentz Spaces. Mediterr. J. Math., **17** (57) (2020). https://doi.org/10.1007/s00009-020-1495-7
5. Bennett, C., Sharpley, R.: Interpolation of operator. Academic Press, Boston (1988)
6. Bergh, J., Löfström, J.: Interpolation Spaces. Springer, Berlin (1976)
7. Boyd, D.W.: Spaces between a pair of reflexive Lebesgue spaces. Proc. Am. Math. Soc. 18, 215–219 (1967)
8. Boyd, D. W.: The Hilbert transform on rearrangement-invariant spaces. Canad. J. Math., **19**, 599 – 616 (1967). DOI: https://doi.org/10.4153/CJM-1967-053-7

9. Boyd, D. W.: Indices of function spaces and their relationship to interpolation, Canad. J. Math., **21**,1245–1254 (1969)
10. Cadilhac, L.: Majorization, Interpolation and noncommutative Khintchine inequalities. Studia Mathematica, **258**(1), 1-26 (2021). 10.4064/sm181217-30-1, hal-02485755
11. Caetano, A., Gogatishvili, A., Opic, B.; Compactness in quasi-Banach function spaces and applications to compact embeddings of Besov-type spaces. Proceedings of the Royal Society of Edinburgh: Section A Mathematics. **146** (5), 905-927 (2016). doi:10.1017/S0308210515000761
12. Dirksen, S.: Noncommutative and Vector-valued Rosenthal Inequalities. PhD thesis, Technische Universiteit Delft. (2011)
13. Dirksen, S.: Noncommutative Boyd interpolation theorems. Transactions of the American Mathematical Society, **367**, 4079-4110 (2015)
14. Doktorskii, R. Ya.: A multiparametric real interpolation method III. Lorentz–Zygmund spaces. Manuscript No. 6070-B88, deposited at VINITI (1988) (Russian)
15. Doktorskii, R. Ya.: A multiparametric real interpolation method IV. Reiteration relations for "limiting" cases $\eta$=0 and $\eta$=1, application to the Fourier series theory. Manuscript No. 4637-B89, deposited at VINITI (1989) (Russian)
16. Doktorskii, R. Ya.: Reiteration relations of the real interpolation method. Soviet. Math. Dokl. **44** (3), 665–669 (1992)
17. Doktorski, L. R. Ya.: An application of limiting interpolation to Fourier series theory. in The Diversity and Beauty of Applied Operator Theory. Operator Theory: Advances and Applications, A. Bottcher, D. Potts, P. Stollmann and D.Wenzel, Eds., **268**, Birkhauser, Cham, Basel (2018)
18. Doktorski, L. R. Ya.: Reiteration formulae for the real interpolation method including $\boldsymbol{\mathcal{L}}$ or $\boldsymbol{\mathcal{R}}$ limiting spaces. J. Funct. Spaces 6685993 (2020). https://doi.org/10.1155/2020/6685993
19. Doktorski, L. R. Ya.: Some reiteration theorems for $\mathcal{R}$, $\mathcal{L}$, $\mathcal{RR}$, $\mathcal{RL}$, $\mathcal{LR}$, and $\mathcal{LL}$ limiting interpolation spaces. J. Funct. Spaces 8513304 (2021). https://doi.org/10.1155/2021/8513304. 8513304.pdf
20. Doktorski, L. R. Ya.: On the spaces introduced by N. K. Karapetyants and B. S. Rubin, and their connection with Lorentz–Zygmund spaces. Bol. Soc. Mat. Mex. **30**(81) (2024). https://doi.org/10.1007/s40590-024-00658-9
21. Doktorski, L. R. Ya.: Every rearrangement-invariant quasi-Banach function space is an interpolation space between two Lorentz spaces. (2026). http://arxiv.org/abs/2511.05096
22. Doktorski, L. R. Ya., Fernández–Martínez, P., Signes, T.: Reiteration theorems for $\mathcal{R}$ and $\mathcal{L}$ spaces with the same parameter. J. Math. Anal. Appl., **508**(1), 125846 (2022). https://doi.org/10.1016/j.jmaa.2021.125846
23. Doktorski, L. R. Ya., Fernández–Martínez, P., Signes, T.: Some examples of equivalent rearrangement-invariant quasi-norms defied via $f^*$ or $f^{**}$. Math. Nachr. **296**(5), 1781–1802 (2023). https://doi.org/10.1002/mana.202200112
24. Doktorski, L. R. Ya., Fernández–Martínez, P., Signes, T.: The $K$-functional and reiteration theorems for left and right spaces, Part I. J. Math. Anal. Appl., **531**(1) 127882 (2024). https://doi.org/10.1016/j.jmaa.2023.127882.
25. Doktorski, L. R. Ya., Fernández–Martínez, P., Signes, T.: The $K$-functional and reiteration theorems for left and right spaces, Part II, J. Math. Anal. Appl., **539**(1) 128482 (2024). https://doi.org/10.1016/j.jmaa.2024.128482

26. Doktorski, L. R. Ya., Fernández–Martínez, P., Signes, T.: Interpolation between grand and small Lorentz–Karamata spaces with the same main parameter. Positivity, **29**(39) (2025). https://doi.org/10.1007/s11117-025-01129-9
27. Edmunds, D. E., Kerman, R., Pick, L.: Optimal Sobolev imbeddings involving rearrangement-invariant quasinorms. J. Funct. Anal., **170**(2), 307–355 (2000)
28. Edmunds, D. E., Opic, B.: Equivalent quasi-norms on Lorentz spaces., Proc. Amer. Math. Soc. **131**, 745–754 (2003)
29. Edmunds, D. E., Opic, B.: Alternative characterisations of Lorentz–Karamata spaces. Czech. Math. J., **58**(2), 517—540 (2008)
30. W. D. Evans and B. Opic, *Real interpolation with logarithmic functors and reiteration*, Canad. J. Math., 52 (5) (2000) 920–960, https://doi.org/10.4153/CJM-2000-039-2.
31. W. D. Evans, B. Opic, and L. Pick, *Real interpolation with logarithmic functors*, J. of lnequal. & Appl., 7 (2) (2002) 187-269. http://dx.doi.org/10.1155/S1025583402000127
32. Fernández–Martínez, P., Signes, T.: Real interpolation with symmetric spaces and slowly varying functions. Quart. J. Math., **63**(1), 133–164 (2012). http://dx.doi.org/10.1093/qmath/haq009
33. Fernández–Martínez, P., Signes, T.: An application of interpolation theory to renorming of Lorentz-Karamata type spaces. Ann. Acad. Sci. Fenn. Math. **39**(1), 97–107 (2014). http://dx.doi.org/10.5186/aasfm.2014.3911
34. Fernández–Martínez, P., Signes, T.: Limit cases of reiteration theorems. Math. Nachr. **288**(1), 25–47 (2015). https://doi.org/10.1002/mana.201300251.
35. Fernández–Martínez, P., Signes, T.: General reiteration theorems for $\mathcal{R}$ and $\mathcal{L}$ classes: Case of left $R$-spaces and right $L$-spaces. J. Math. Anal. Appl., **494**, 124649 (2021). https://doi.org/10.1016/j.jmaa.2020.124649.
36. Fernández–Martínez, P., Signes, T.: General reiteration theorems for $\mathcal{R}$ and $\mathcal{L}$ classes: mixed interpolation of $\mathcal{R}$ and $\mathcal{L}$-spaces. Positivity, **26**(47) (2022). https://doi.org/10.1007/s11117-022-00888-z
37. Fiorenza, A.: Duality and reflexivity in grand Lebesgue spaces. Collect. Math. **51**, 131–148 (2000)
38. Fiorenza, A., Formica, M. R., Gogatishvili, A., Kopaliani, T., Rakotoson, J. M.: Characterization of interpolation between grand, small or classical Lebesgue spaces. Nonlinear Analysis, **177**, 422–453 (2018). https://doi.org/10.1016/j.na.2017.09.005
39. Fiorenza, A., Karadzhov, G.E.: Grand and small Lebesgue spaces and their analogs. Z. Anal. Anwend. **23**(4), 657–681 (2004). https://doi.org/10.4171/zaa/1215
40. Formica, M. R., Giova, R.: Boyd indices in generalized grand Lebesgue spaces and applications. Mediterr. J. Math. **12**, 987–995 (2015) https://doi.org/10.1007/s00009-014-0439-5.
41. Di Fratta, G., Fiorenza, A.: A direct approach to the duality of grand and small Lebesgue spaces. Nonlinear Analysis **70**, 2582–2592 (2009). https://doi.org/10.1016/j.na.2008.03.044
42. Gogatishvili, A., Ünver, T.: Weighted inequalities involving two Hardy operators. Anal. Math. Phys., **15**(99) (2025). https://doi.org/10.1007/s13324-025-01101-6
43. Gogatishvili, A., Mustafayev, R., Ünver, T.: Embeddings between weighted Copson and Cesàro function spaces, Czech. Math. J., **67**(142), 1105–1132 (2017). https://doi.org/10.21136/CMJ.2017.0424-16

44. Gogatishvili, A., Opic, B., Trebels, W.: Limiting reiteration for real interpolation with slowly varying functions. Math. Nachr. **278**, 86–107 (2005). doi.org/10.1002/mana.200310228
45. Hudzik, H., Maligranda, L.: An interpolation theorem in symmetric function F-spaces. Proceedings of the American Mathematical Society, **110**(1), 89-96 (1990)
46. Iwaniec, T., Sbordone, C.: On the integrability of the Jacobian under minimal hypotheses. Arch. Rational Mech. Anal. **119**, 129–143 (1992)
47. Kamińska, A., Raynaud, Y.: Isomorphic copies in the lattice $E$ and its symmetrization $E(*)$ with applications to Orlicz-Lorentz spaces. J. Funct. Anal. **257**(1), 271-331 (2009)
48. Karadzhov, G. E.: Boyd indices for quasi-normed rearrangement invariant spaces. Journal of Prime Research in Mathematics **8**, 36-44 (2012)
49. Kolwicz, P., Leśnik, K., Maligranda, L.: Pointwise products of some Banach function spaces and factorization. J. Funct. Anal. **266**(2), 616-659 (2014)
50. Kolwicz, P., Leśnik, K., Maligranda, L.: Symmetrization, factorization and arithmetic of quasi-Banach function spaces. J. Math. Anal. Appl., **470**(2), 1136-1166 (2019). https://doi.org/10.1016/j.jmaa.2018.10.054
51. Krein, S. G., Petunin, Ju. I., Semenov, E. M.: Interpolation of linear operators, Translations of Math. Monographs, **54**, Amer. Math. Soc., Providence R.I., (1982)
52. Kufner, A., Maligranda, L., Persson, L.-E.: The Hardy inequality: about its history and some related results. Vydavatelský servis, Praha (2007)
53. Lorentz, G. G.: On the theory of spaces $\Lambda$. Pacific J. Math. **1**, 411-429 (1951)
54. Lorist, E., Nieraeth, Z.: Banach function spaces done right. Indag. Math. (N.S.), **35**(2), 247-268 (2024). https://doi.org/10.1016/j.indag.2023.11.004
55. Luxemburg, W. A. J.: Banach function spaces. Thesis, Delft Technical Univ. (1955)
56. Luxemburg, W. A. J.: Rearrangement-invariant Banach function spaces. In Proc. Sympos. in Analysis, volume 10 of Queen's Papers in Pure and Appl. Math., 83–144. Queen's University (1967)
57. Maligranda, L.: Generalized Hardy inequalities in rearrangement invariant spaces, J. Math. Pures Appl. **59**, 405–415 (1980)
58. Maligranda, L.: Indices and interpolation. Dissertationes Math. (Rozprawy Mat.), **234**(49) (1985)
59. Montgomery-Smith, S. J.: The Hardy operator and Boyd indices. In Interaction between functional analysis, harmonic analysis, and probability (Columbia, MO, 1994), volume 175 of Lecture Notes in Pure and Appl. Math., 359–364. Dekker, New York.
60. Musilová, A., Nekvinda, A., Peša, D., Turčinová, H.: On the properties of rearrangement-invariant quasi-Banach function spaces. (2024). https://doi.org/10.48550/arXiv.2404.00707
61. Nekvinda, A. Peša, D.: On the properties of quasi-Banach function spaces, J. Geom. Anal. **34**(231) (2024). https://doi.org/10.1007/s12220-024-01673-y
62. Opic, B.: New characterizations of Lorentz spaces. Proc. Royal Soc. Edinburgh **133A**, 439–448 (2003)
63. Opic, B.: On equivalent quasi-norms on Lorentz spaces. In Function Spaces, Differential Operators and Nonlinear Analysis. The Hans Triebel Aniversary Volume (D. Haroske, T. Runst, H.-J. Schmeisser, eds.). Birkhäuser Verlag, Basel/Switzerland, 415–426 2003

64. Opic, B., Pick, L,: On generalized Lorentz-Zygmund spaces. Math. Ineq. Appl. **2**, 391–467 (1999). https://doi.org/10.7153/mia-02-35
65. Peša, D.: Lorentz–Karamata spaces. (2023). https://arxiv.org/abs/2006.14455v2
66. Peša, D.: On the smoothness of slowly varying functions. Proceedings of the Edinburgh Mathematical Society, **67**(3), 876-891 (2024). https://doi.org/10.1017/S0013091524000348
67. Sagher, Y.: Interpolation of *r*-Banach spaces. Stud. Math., **41**(1), 45-70 (1972)
68. Sawano, Y., Ho, K.-P., Yang, D., Yang, S.: Hardy spaces for ball quasi-Banach function spaces. Dissertat. Math. **525**(102) (2017)
69. Triebel, H.: Interpolation theory, function spaces, differential operators. North-Holland, Amsterdam (1978)
70. Yoshinaga, K.: Lorentz spaces and the Calderón–Zygmund theorem. Bull. Kyushu Inst. Tech. **16**, 1-38 (1969)